\newcommand{\xbD}{\Delta}

\newcommand{\xbG}{\Gamma}

\newcommand{\xbP}{\Pi}

\newcommand{\xba}{\alpha}
\newcommand{\xbb}{\beta}

\newcommand{\xbd}{\delta}
\newcommand{\xbe}{\in}
\newcommand{\xbf}{\phi}
\newcommand{\xbg}{\gamma}

\newcommand{\xbk}{\kappa}
\newcommand{\xbl}{\lambda}
\newcommand{\xbm}{\mu}

\newcommand{\xbo}{\omega}

\newcommand{\xbq}{\psi}

\newcommand{\xbs}{\sigma}

\newcommand{\xCK}{\times}

\newcommand{\xCN}{\neg}
\newcommand{\xCQ}{\emptyset}

\newcommand{\xCf}{\hspace{0.1em}}

\newcommand{\xcA}{\forall}

\newcommand{\xcC}{\not\subseteq}

\newcommand{\xcE}{\exists}

\newcommand{\xcH}{\not\Rightarrow}
\newcommand{\xcI}{\not\Leftarrow}

\newcommand{\xcN}{\hspace{0.2em}\not\sim\hspace{-0.9em}\mid\hspace{0.8em}}
\newcommand{\xcO}{\bigvee}
\newcommand{\xcP}{\not\rightarrow}

\newcommand{\xcS}{\bigcap}
\newcommand{\xcT}{\bot}
\newcommand{\xcU}{\bigwedge}
\newcommand{\xcV}{\bigcup}

\newcommand{\xca}{\infty}
\newcommand{\xcb}{\subset}
\newcommand{\xcc}{\subseteq}
\newcommand{\xcd}{\supseteq}
\newcommand{\xce}{\not\in}
\newcommand{\xcf}{\supset}
\newcommand{\xcg}{\geq}
\newcommand{\xch}{\Rightarrow}
\newcommand{\xci}{\Leftarrow}
\newcommand{\xcj}{\Leftrightarrow}
\newcommand{\xck}{\leq}
\newcommand{\xcl}{\vdash}
\newcommand{\xcm}{\models}
\newcommand{\xcn}{\hspace{0.2em}\sim\hspace{-0.9em}\mid\hspace{0.58em}}

\newcommand{\xco}{\vee}
\newcommand{\xcp}{\rightarrow}

\newcommand{\xcr}{\leftrightarrow}
\newcommand{\xcs}{\cap}
\newcommand{\xcu}{\wedge}
\newcommand{\xcv}{\cup}

\newcommand{\xcz}{\Box}

\newcommand{\xDH}{\item }

\newcommand{\xdC}{\mbox{\boldmath$C$}}
\newcommand{\xdD}{\mbox{\boldmath$D$}}

\newcommand{\xdl}{{\cal L}}
\newcommand{\xdm}{{\cal M}}

\newcommand{\xdp}{{\cal P}}

\newcommand{\xdu}{{\cal U}}

\newcommand{\xdx}{{\cal X}}
\newcommand{\xdy}{{\cal Y}}

\newcommand{\xEH}{ & }
\newcommand{\xEI}{\begin{itemize}}
\newcommand{\xEJ}{\end{itemize}}
\newcommand{\xEP}{ \\ }

\newcommand{\xEd}{\neq}
\newcommand{\xEh}{\begin{enumerate}}
\newcommand{\xEj}{\end{enumerate}}

\newcommand{\xeb}{\prec}

\newcommand{\xex}{\lceil}

\newcommand{\xFO}{\parallel}

\newcommand{\Xl}{\ldots}

\newcommand{\ol}{\overline}
\newcommand{\ul}{\underline}

\newcommand{\xssc}{\scriptsize}

\newcommand{\bl}{\begin{lemma} \rm}
\newcommand{\el}{\end{lemma}}
\newcommand{\br}{\begin{remark} \rm}
\newcommand{\er}{\end{remark}}
\newcommand{\be}{\begin{example} \rm}
\newcommand{\ee}{\end{example}}
\newcommand{\bco}{\begin{corollary} \rm}
\newcommand{\eco}{\end{corollary}}
\newcommand{\bc}{\begin{claim} \rm}
\newcommand{\ec}{\end{claim}}
\newcommand{\bfa}{\begin{fact} \rm}
\newcommand{\efa}{\end{fact}}
\newcommand{\bp}{\begin{proposition} \rm}
\newcommand{\ep}{\end{proposition}}
\newcommand{\bd}{\begin{definition} \rm}
\newcommand{\ed}{\end{definition}}
\newcommand{\bcs}{\begin{construction} \rm}
\newcommand{\ecs}{\end{construction}}
\newcommand{\bcd}{\begin{condition} \rm}
\newcommand{\ecd}{\end{condition}}
\newcommand{\bt}{\begin{theorem} \rm}
\newcommand{\et}{\end{theorem}}
\newcommand{\bn}{\begin{notation} \rm}
\newcommand{\en}{\end{notation}}
\newcommand{\bfi}{\begin{bild} \rm}
\newcommand{\efi}{\end{bild}}
\newcommand{\bsta}{\begin{statement} \rm}
\newcommand{\esta}{\end{statement}}
\newcommand{\bcom}{\begin{comment} \rm}
\newcommand{\ecom}{\end{comment}}
\newcommand{\bdia}{\begin{diagram} \rm}
\newcommand{\edia}{\end{diagram}}

\newcommand{\bfc}{\begin{figure}[htb] \begin{center}}
\newcommand{\efc}{\end{center} \end{figure}}

\documentclass{article}

\usepackage{amssymb,epic,eepic}

\title{
Cumulativity without closure of the domain under finite unions
}

\author{Dov M Gabbay
\thanks{
Dov.Gabbay@kcl.ac.uk, www.dcs.kcl.ac.uk/staff/dg
} \\
King's College, London
\thanks{
Department of Computer Science, King's College London, Strand,
London WC2R 2LS, UK
} \\ \\
Karl Schlechta
\thanks{
ks@cmi.univ-mrs.fr, karl.schlechta@web.de, http://www.cmi.univ-mrs.fr/ $\sim$ ks
} \\
Laboratoire d'Informatique Fondamentale de Marseille
\thanks{
UMR 6166, CNRS and Universit\'{e} de Provence,
Address: CMI, 39, rue Joliot-Curie, F-13453 Marseille Cedex 13, France
}
}

\begin{document}

\newtheorem{lemma}{Lemma}[section]
\newtheorem{theorem}[lemma]{Theorem}
\newtheorem{proposition}[lemma]{Proposition}
\newtheorem{corollary}[lemma]{Corollary}
\newtheorem{claim}[lemma]{Claim}
\newtheorem{fact}[lemma]{Fact}
\newtheorem{remark}[lemma]{Remark}
\newtheorem{definition}{Definition}[section]
\newtheorem{construction}{Construction}[section]
\newtheorem{condition}{Condition}[section]
\newtheorem{example}{Example}[section]
\newtheorem{notation}{Notation}[section]
\newtheorem{bild}{Figure}[section]
\newtheorem{comment}{Comment}[section]
\newtheorem{statement}{Statement}[section]
\newtheorem{diagram}{Diagram}[section]

\maketitle

\renewcommand{\labelenumi}
  {(\arabic{enumi})}
\renewcommand{\labelenumii}
  {(\arabic{enumi}.\arabic{enumii})}
\renewcommand{\labelenumiii}
  {(\arabic{enumi}.\arabic{enumii}.\arabic{enumiii})}
\renewcommand{\labelenumiv}
  {(\arabic{enumi}.\arabic{enumii}.\arabic{enumiii}.\arabic{enumiv})}

\setcounter{secnumdepth}{3}
\setcounter{tocdepth}{3}

\begin{abstract}

For nonmonotonic logics, Cumulativity is an important logical rule.
We show here that Cumulativity fans out into
an infinity of different conditions, if the domain is not closed under
finite unions.

\end{abstract}

\tableofcontents

%
%
%
\section{Introduction}
\subsection{
Motivation and history of Cumulativity
}

Cumulativity is one of the
most important properties a nonmonotonic logic can have, and was
recognised as
such quite long ago, see  \cite{Gab85}. It says that
adding results proved already to the axioms will not change the set of
provable results. More precisely, if $ \xcn $ is some consequence
relation,
then cumulativity says:

If $ \xbf \xcn \xbq,$ then: $ \xbf \xcn \xbq ' $ iff $ \xbf \xcu \xbq
\xcn \xbq '.$

(The other axiom in Gabbay's system is $ \xba \xcn \xba.$ - We modified
both
slightly here, the precise definition, and other axioms as well as
variants are
given in Definition \ref{Definition Log-Cond}.)

In classical logic, this is trivial, in nonmonotonic
logics, this is a non-trivial property.

Whereas the approach of  \cite{Gab85} was an abstract consideration
of
desirable properties, preferential structures were introduced as
abstractions
of Circumscription independently in  \cite{Sho87b} and  \cite{BS85}. Roughly,
a preferential model is a possible worlds structure with a binary
relation $ \xeb,$ which expresses that one world is more ``normal'' than
the
other. A precise definition is given below in Definition \ref{Definition
Pref-Str}.

Both approaches were connected in  \cite{KLM90}, where a
representation
theorem was proved, showing that Cumulativity corresponds to
the relational property of ``smoothness'', which says that every possible
world is
either minimal itself, or there is a smaller one, which is minimal. This
property can be violated by non-transitive relations, or infinite
descending
chains.

Preferential semantics generate, however, a richer logic than the one
considered in  \cite{Gab85}, and, to the authors' knowledge, a
precise
semantics for this system is given only in  \cite{GS08b} by the
present authors.

Cumulativity can also be seen as allowing certain manipulations
with ``small'' sets. If $ \xbf \xcn \xbq $ expresses that the number of
exceptions,
i.e. when $ \xbf \xcu \xCN \xbq $ holds, is small, then Cumulativity
essentiall says that
the number of cases where $ \xbf \xcu \xbq \xcu \xCN \xbq ' $ holds is
still small in the
set of cases where $ \xbf \xcu \xbq $ holds, or, if A and $A' $ are small
subsets of $B,$
then $A' $ will still be small in $B-A$ (assume for simplicity that A and
$A' $
are disjoint). See Diagram \ref{Diagram CumSmall}. A systematic
investigation of the connections between ``size'' and nonmonotonic logic
can be found in  \cite{GS08c}.

\vspace{10mm}

\begin{diagram}

\label{Diagram CumSmall}
\index{Diagram CumSmall}

\centering
\setlength{\unitlength}{1mm}
{\renewcommand{\dashlinestretch}{30}
\begin{picture}(100,100)(0,0)
\put(50,50){\circle{80}}
\path(24.3,80.6)(75.7,80.6)
\path(84.64,70)(84.64,30)

\put(50,50){{\xssc $B$}}
\put(50,85){{\xssc $A$}}
\put(85,50){{\xssc $A'$}}

\put(30,5) {{\rm\bf Cumulativity}}

\end{picture}
}
\end{diagram}

\vspace{4mm}

Quite often, the domain on which a logic operates is closed under finite
union: if $M( \xbf )$ is the set of $ \xbf -$models, then $M( \xbf ) \xcv
M( \xbf ' )=M( \xbf \xco \xbf ' ).$
This need not always be the case, for instance, not for certain sequent
calculi - see Section \ref{Section Sequent-Calculi} below - which do not
have
an ``or'' in the language.

(Yet, even if the logic allows such closure, i.e. we have ``or'' in the
language,
we might simply be unable to observe the result of the closure,
we might well know the consequences of $ \xbf,$ those of $ \xbf ',$ but
be unable to
observe the consequences of $ \xbf \xco \xbf '.$ We may think here of an
experiment in,
e.g., physics,
where the input is a formula, and the output is an observation, again a
formula. Our experimental setup may just not allow to observe the outcome
of $ \xbf \xco \xbf ' $ - though the language allows to pose the question
$'' \xbf \xco \xbq?''.)$

Quite surprisingly (for the authors at least), absence of this simple
closure
condition has drastic consequences on cumulativity: the simple short
condition
breaks up into an infinity of non-equivalent conditions. This is the
main result of this paper, and shown
by Example \ref{Example Inf-Cum-Alpha}.
For illustration, we now give the first three cases of one set of
conditions, the full definition is given
in Definition \ref{Definition Cum-Alpha}.

We switch to semantics, and give a corresponding
condition for model sets. Often, it is good policy to separate the
semantical from the proof theoretical conditions, as both have their
own problems, and it is easier to treat them separately. More
connections are given in Definition \ref{Definition Log-Cond}. $ \xbm (X)$
will now be the
set of minimal models of $X.$

 \xEh
 \xDH
$( \xbm Cum0)$ $ \xbm (X_{0}) \xcc U$ $ \xcp $ $X_{0} \xcs \xbm (U) \xcc
\xbm (X_{0}),$

 \xDH
$( \xbm Cum1)$ $ \xbm (X_{0}) \xcc U,$ $ \xbm (X_{1}) \xcc U \xcv X_{0}$ $
\xcp $ $X_{0} \xcs X_{1} \xcs \xbm (U) \xcc \xbm (X_{1}),$

 \xDH
$( \xbm Cum2)$ $ \xbm (X_{0}) \xcc U,$ $ \xbm (X_{1}) \xcc U \xcv X_{0},$
$ \xbm (X_{2}) \xcc U \xcv X_{0} \xcv X_{1}$ $ \xcp $ $X_{0} \xcs X_{1}
\xcs X_{2} \xcs \xbm (U) \xcc \xbm (X_{2}).$

 \xEj

The semantical version of traditional cumulativity is

$( \xbm CUM)$ $ \xbm (X) \xcc Y \xcc X$ $ \xch $ $ \xbm (X)= \xbm (Y).$

It is easy to show that $( \xbm CUM)$ and closure under finite unions
imply
$( \xbm Cumi),$ but one has to be careful about the prerequisites used, so
we refer the reader to Fact \ref{Fact Cum-Alpha} (5.3).

Let $T$ etc. be sets of formulas, and assume them to be closed under
classical deduction, let $M(T)$ be the set of classical models of $T.$
Assume further that the model choice function $ \xbm $ has the property
that
$ \xbm (M(T))=M(T' )$ for some $T' $ (i.e. $ \xbm $ is definability
preserving
in the sense of Definition \ref{Definition Log-Base}), and write
$ \ol{ \ol{T} }$ for $\{ \xbf:T \xcn \xbf \}.$

Then these semantical conditions for a model choice function $ \xbm $
translate
as follows into logic:

 \xEh
 \xDH
(Cum0) $T \xcc \ol{ \ol{T_{0}} }$ $ \xch $ $ \ol{ \ol{T_{0}} } \xcc T_{0}
\xcv \ol{ \ol{T} },$

 \xDH
(Cum1) $T \xcc \ol{ \ol{T_{0}} },$ $T \xcs T_{0} \xcc \ol{ \ol{T_{1}} }$ $
\xch $ $ \ol{ \ol{T_{1}} } \xcc T_{0} \xcv T_{1} \xcv \ol{ \ol{T} }$

 \xDH
(Cum2) $T \xcc \ol{ \ol{T_{0}} },$ $T \xcs T_{0} \xcc \ol{ \ol{T_{1}} },$
$T \xcs T_{0} \xcs T_{1} \xcc \ol{ \ol{T_{2}} }$ $ \xch $ $ \ol{
\ol{T_{2}} } \xcc T_{0} \xcv T_{1} \xcv T_{2} \xcv \ol{
\ol{T} }$

 \xEj

Cumulativity for formula sets is now (slightly simplified)

$ \xCf (CUM)$ $T \xcc T' \xcc \ol{ \ol{T} }$ $ \xch $ $ \ol{ \ol{T} }=
\ol{ \ol{T' } }.$

Seen conversely, closure conditions of the
domain - or, better, lack of them - thus reveal themselves as powerful
tools to
separate logical rules.
\subsection{Organisation of the paper}

We first give the general framework we work in:
We mention the logical properties often considered in the framework of
nonmonotonic, and in particular preferential logics, as well as their
algebraic
counterparts. The reader will
find here more than he really needs to understand the main result, but it
will
help him to put things into perspective. We also show (without proof)
connections and differences, which can be quite subtle - see
[GS08c] for details and proofs.

We then present two sequent calculi, one due to Lehmann, the other to
Arieli and Avron. The first one demonstrated the authors that domain
closure
questions can be real problems, and were not just due to an incapacity to
find a proof with weaker prerequisites. The crucial point in both is that
there is no ``or'' in the language, so the domain is not closed under finite
unions.

The main result of the paper, shown in
Section \ref{Section Cumulativity-without-union},
Example \ref{Example Inf-Cum-Alpha}
demonstrates that the problem is quite serious, as it splits the otherwise
simple condition of Cumulativity into an infinity of conditions, which all
collaps to one in the presence of closure of the domain under finite
unions.
We conclude by giving positive results for the different conditions of
cumulativity - see Fact \ref{Fact Cum-Alpha}.
\section{
Basic definitions for preferential structures and related logics
}
\label{Section BasicDefinitions}
\index{Definition Alg-Base}

\bd

$\hspace{0.01em}$


\label{Definition Alg-Base}

We use $ \xdp $ to denote the power set operator,
$ \xbP \{X_{i}:i \xbe I\}$ $:=$ $\{g:$ $g:I \xcp \xcV \{X_{i}:i \xbe I\},$
$ \xcA i \xbe I.g(i) \xbe X_{i}\}$ is the general cartesian
product, $card(X)$ shall denote the cardinality of $X,$ and $V$ the
set-theoretic
universe we work in - the class of all sets. Given a set of pairs $ \xdx
,$ and a
set $X,$ we denote by $ \xdx \xex X:=\{<x,i> \xbe \xdx:x \xbe X\}.$ When
the context is clear, we
will sometime simply write $X$ for $ \xdx \xex X.$

$A \xcc B$ will denote that $ \xCf A$ is a subset of $B$ or equal to $B,$
and $A \xcb B$ that $ \xCf A$ is
a proper subset of $B,$ likewise for $A \xcd B$ and $A \xcf B.$

Given some fixed set $U$ we work in, and $X \xcc U,$ then $ \xdC (X):=U-X$
.

If $ \xdy \xcc \xdp (X)$ for some
$X,$ we say that $ \xdy $ satisfies

$( \xcs )$ iff it is closed under finite intersections,

$( \xcS )$ iff it is closed under arbitrary intersections,

$( \xcv )$ iff it is closed under finite unions,

$( \xcV )$ iff it is closed under arbitrary unions,

$( \xdC )$ iff it is closed under complementation.

We will sometimes write $A=B \xFO C$ for: $A=B,$ or $A=C,$ or $A=B \xcv
C.$

We make ample and tacit use of the Axiom of Choice.
\index{Definition Log-Base}

\ed

\bd

$\hspace{0.01em}$


\label{Definition Log-Base}

We work here in a classical propositional language $ \xdl,$ a theory $T$
will be an
arbitrary set of formulas. Formulas will often be named $ \xbf,$ $ \xbq
,$ etc., theories
$T,$ $S,$ etc.

$v( \xdl )$ will be the set of propositional variables of $ \xdl.$

$M_{ \xdl }$ will be the set of (classical) models of $ \xdl,$ $M(T)$ or
$M_{T}$
is the set of models of $T,$ likewise $M( \xbf )$ for a formula $ \xbf.$

$ \xdD_{ \xdl }:=\{M(T):$ $T$ a theory in $ \xdl \},$ the set of definable
model sets.

Note that, in classical propositional logic, $ \xCQ,M_{ \xdl } \xbe
\xdD_{ \xdl },$ $ \xdD_{ \xdl }$ contains
singletons, is closed under arbitrary intersections and finite unions.

An operation $f: \xdy \xcp \xdp (M_{ \xdl })$ for $ \xdy \xcc \xdp (M_{
\xdl })$ is called definability
preserving, $ \xCf (dp)$ or $( \xbm dp)$ in short, iff for all $X \xbe
\xdD_{ \xdl } \xcs \xdy $ $f(X) \xbe \xdD_{ \xdl }.$

We will also use $( \xbm dp)$ for binary functions $f: \xdy \xCK \xdy \xcp
\xdp (M_{ \xdl })$ - as needed
for theory revision - with the obvious meaning.

$ \xcl $ will be classical derivability, and

$ \ol{T}:=\{ \xbf:T \xcl \xbf \},$ the closure of $T$ under $ \xcl.$

$Con(.)$ will stand for classical consistency, so $Con( \xbf )$ will mean
that
$ \xbf $ is clasical consistent, likewise for $Con(T).$ $Con(T,T' )$ will
stand for
$Con(T \xcv T' ),$ etc.

Given a consequence relation $ \xcn,$ we define

$ \ol{ \ol{T} }:=\{ \xbf:T \xcn \xbf \}.$

(There is no fear of confusion with $ \ol{T},$ as it just is not useful to
close
twice under classical logic.)

$T \xco T':=\{ \xbf \xco \xbf ': \xbf \xbe T, \xbf ' \xbe T' \}.$

If $X \xcc M_{ \xdl },$ then $Th(X):=\{ \xbf:X \xcm \xbf \},$ likewise
for $Th(m),$ $m \xbe M_{ \xdl }.$
\index{Definition Log-Cond}

\ed

\bd

$\hspace{0.01em}$


\label{Definition Log-Cond}

We introduce here formally a list of properties of set functions on the
algebraic side, and their corresponding logical rules on the other side.

Recall that $ \ol{T}:=\{ \xbf:T \xcl \xbf \},$ $ \ol{ \ol{T} }:=\{ \xbf
:T \xcn \xbf \},$
where $ \xcl $ is classical consequence, and $ \xcn $ any other
consequence.

We show, wherever adequate, in parallel the formula version
in the left column, the theory version
in the middle column, and the semantical or algebraic
counterpart in the
right column. The algebraic counterpart gives conditions for a
function $f:\xdy\xcp\xdp (U)$, where $U$ is some set, and
$\xdy\xcc\xdp (U)$.

When the formula version is not commonly used, we omit it,
as we normally work only with the theory version.

Intuitively, $A$ and $B$ in the right hand side column stand for
$M(\xbf)$ for some formula $\xbf$, whereas $X$, $Y$ stand for
$M(T)$ for some theory $T$.

{\footnotesize

\begin{tabular}{|c|c|c|}

\hline

\multicolumn{3}{|c|}{Basics} \xEP

\hline

$(AND)$
\xEH
$(AND)$
\xEH
Closure under
\xEP

$ \xbf \xcn \xbq,  \xbf \xcn \xbq '   \xch $
\xEH
$ T \xcn \xbq, T \xcn \xbq '   \xch $
\xEH
finite
\xEP

$ \xbf \xcn \xbq \xcu \xbq ' $
\xEH
$ T \xcn \xbq \xcu \xbq ' $
\xEH
intersection
\xEP

\hline

$(OR)$ \xEH $(OR)$ \xEH $( \xbm OR)$ \xEP

$ \xbf \xcn \xbq,  \xbf ' \xcn \xbq   \xch $ \xEH
$ \ol{\ol{T}} \xcs \ol{\ol{T'}} \xcc \ol{\ol{T \xco T'}} $ \xEH
$f(X \xcv Y) \xcc f(X) \xcv f(Y)$
\xEP

$ \xbf \xco \xbf ' \xcn \xbq $ \xEH
\xEH
\xEP

\hline

$(wOR)$
\xEH
$(wOR)$
\xEH
$( \xbm wOR)$
\xEP

$ \xbf \xcn \xbq,$ $ \xbf ' \xcl \xbq $ $ \xch $
\xEH
$ \ol{ \ol{T} } \xcs \ol{T' }$ $ \xcc $ $ \ol{ \ol{T \xco T' } }$
\xEH
$f(X \xcv Y) \xcc f(X) \xcv Y$
\xEP

$ \xbf \xco \xbf ' \xcn \xbq $
\xEH
\xEH
\xEP

\hline

$(disjOR)$
\xEH
$(disjOR)$
\xEH
$( \xbm disjOR)$
\xEP

$ \xbf \xcl \xCN \xbf ',$ $ \xbf \xcn \xbq,$
\xEH
$\xCN Con(T \xcv T') \xch$
\xEH
$X \xcs Y= \xCQ $ $ \xch $
\xEP

$ \xbf ' \xcn \xbq $ $ \xch $ $ \xbf \xco \xbf ' \xcn \xbq $
\xEH
$ \ol{ \ol{T} } \xcs \ol{ \ol{T' } } \xcc \ol{ \ol{T \xco T' } }$
\xEH
$f(X \xcv Y) \xcc f(X) \xcv f(Y)$
\xEP

\hline

$(LLE)$
\xEH
$(LLE)$
\xEH
\xEP

Left Logical Equivalence
\xEH
\xEH
\xEP

$ \xcl \xbf \xcr \xbf ',  \xbf \xcn \xbq   \xch $
\xEH
$ \ol{T}= \ol{T' }  \xch   \ol{\ol{T}} = \ol{\ol{T'}}$
\xEH
trivially true
\xEP

$ \xbf ' \xcn \xbq $ \xEH \xEH \xEP

\hline

$(RW)$ Right Weakening
\xEH
$(RW)$
\xEH
upward closure
\xEP

$ \xbf \xcn \xbq,  \xcl \xbq \xcp \xbq '   \xch $
\xEH
$ T \xcn \xbq,  \xcl \xbq \xcp \xbq '   \xch $
\xEH
\xEP

$ \xbf \xcn \xbq ' $
\xEH
$T \xcn \xbq ' $
\xEH
\xEP

\hline

$(CCL)$ Classical Closure \xEH $(CCL)$ \xEH \xEP

\xEH
$ \ol{ \ol{T} }$ is classically
\xEH
trivially true
\xEP

\xEH closed \xEH \xEP

\hline

$(SC)$ Supraclassicality \xEH $(SC)$ \xEH $( \xbm \xcc )$ \xEP

$ \xbf \xcl \xbq $ $ \xch $ $ \xbf \xcn \xbq $ \xEH $ \ol{T} \xcc \ol{
\ol{T} }$ \xEH $f(X) \xcc X$ \xEP

\cline{1-1}

$(REF)$ Reflexivity \xEH \xEH \xEP
$ \xbD,\xba \xcn \xba $ \xEH \xEH \xEP

\hline

$(CP)$ \xEH $(CP)$ \xEH $( \xbm \xCQ )$ \xEP

Consistency Preservation \xEH \xEH \xEP

$ \xbf \xcn \xcT $ $ \xch $ $ \xbf \xcl \xcT $ \xEH $T \xcn \xcT $ $ \xch
$ $T \xcl \xcT $ \xEH $f(X)= \xCQ $ $ \xch $ $X= \xCQ $ \xEP

\hline

\xEH
\xEH $( \xbm \xCQ fin)$
\xEP

\xEH
\xEH $X \xEd \xCQ $ $ \xch $ $f(X) \xEd \xCQ $
\xEP

\xEH
\xEH for finite $X$
\xEP

\hline

\xEH $(PR)$ \xEH $( \xbm PR)$ \xEP

$ \ol{ \ol{ \xbf \xcu \xbf ' } }$ $ \xcc $ $ \ol{ \ol{ \ol{ \xbf } } \xcv
\{ \xbf ' \}}$ \xEH
$ \ol{ \ol{T \xcv T' } }$ $ \xcc $ $ \ol{ \ol{ \ol{T} } \xcv T' }$ \xEH
$X \xcc Y$ $ \xch $
\xEP

\xEH \xEH $f(Y) \xcs X \xcc f(X)$
\xEP

\cline{3-3}

\xEH
\xEH
$(\xbm PR ')$
\xEP

\xEH
\xEH
$f(X) \xcs Y \xcc f(X \xcs Y)$
\xEP

\hline

$(CUT)$ \xEH $(CUT)$ \xEH $ (\xbm CUT) $ \xEP
$ \xbD \xcn \xba; \xbD, \xba \xcn \xbb \xch $ \xEH
$T \xcc \ol{T' } \xcc \ol{ \ol{T} }  \xch $ \xEH
$f(X) \xcc Y \xcc X  \xch $ \xEP
$ \xbD \xcn \xbb $ \xEH
$ \ol{ \ol{T'} } \xcc \ol{ \ol{T} }$ \xEH
$f(X) \xcc f(Y)$
\xEP

\hline

\end{tabular}

}

{\footnotesize

\begin{tabular}{|c|c|c|}

\hline

\multicolumn{3}{|c|}{Cumulativity} \xEP

\hline

$(CM)$ Cautious Monotony \xEH $(CM)$ \xEH $ (\xbm CM) $ \xEP

$ \xbf \xcn \xbq,  \xbf \xcn \xbq '   \xch $ \xEH
$T \xcc \ol{T' } \xcc \ol{ \ol{T} }  \xch $ \xEH
$f(X) \xcc Y \xcc X  \xch $
\xEP

$ \xbf \xcu \xbq \xcn \xbq ' $ \xEH
$ \ol{ \ol{T} } \xcc \ol{ \ol{T' } }$ \xEH
$f(Y) \xcc f(X)$
\xEP

\cline{1-1}

\cline{3-3}

or $(ResM)$ Restricted Monotony \xEH \xEH $(\xbm ResM)$ \xEP
$ \xbD \xcn \xba, \xbb \xch \xbD,\xba \xcn \xbb $ \xEH \xEH
$ f(X) \xcc A \xcs B \xch f(X \xcs A) \xcc B $ \xEP

\hline

$(CUM)$ Cumulativity \xEH $(CUM)$ \xEH $( \xbm CUM)$ \xEP

$ \xbf \xcn \xbq   \xch $ \xEH
$T \xcc \ol{T' } \xcc \ol{ \ol{T} }  \xch $ \xEH
$f(X) \xcc Y \xcc X  \xch $
\xEP

$( \xbf \xcn \xbq '   \xcj   \xbf \xcu \xbq \xcn \xbq ' )$ \xEH
$ \ol{ \ol{T} }= \ol{ \ol{T' } }$ \xEH
$f(Y)=f(X)$ \xEP

\hline

\xEH
$ (\xcc \xcd) $
\xEH
$ (\xbm \xcc \xcd) $
\xEP
\xEH
$T \xcc \ol{\ol{T'}}, T' \xcc \ol{\ol{T}} \xch $
\xEH
$ f(X) \xcc Y, f(Y) \xcc X \xch $
\xEP
\xEH
$ \ol{\ol{T'}} = \ol{\ol{T}}$
\xEH
$ f(X)=f(Y) $
\xEP

\hline

\multicolumn{3}{|c|}{Rationality} \xEP

\hline

$(RatM)$ Rational Monotony \xEH $(RatM)$ \xEH $( \xbm RatM)$ \xEP

$ \xbf \xcn \xbq,  \xbf \xcN \xCN \xbq '   \xch $ \xEH
$Con(T \xcv \ol{\ol{T'}})$, $T \xcl T'$ $ \xch $ \xEH
$X \xcc Y, X \xcs f(Y) \xEd \xCQ   \xch $
\xEP

$ \xbf \xcu \xbq ' \xcn \xbq $ \xEH
$ \ol{\ol{T}} \xcd \ol{\ol{\ol{T'}} \xcv T} $ \xEH
$f(X) \xcc f(Y) \xcs X$ \xEP

\hline

\xEH $(RatM=)$ \xEH $( \xbm =)$ \xEP

\xEH
$Con(T \xcv \ol{\ol{T'}})$, $T \xcl T'$ $ \xch $ \xEH
$X \xcc Y, X \xcs f(Y) \xEd \xCQ   \xch $
\xEP

\xEH
$ \ol{\ol{T}} = \ol{\ol{\ol{T'}} \xcv T} $ \xEH
$f(X) = f(Y) \xcs X$ \xEP

\hline

\xEH
$(Log=' )$
\xEH $( \xbm =' )$
\xEP

\xEH
$Con( \ol{ \ol{T' } } \xcv T)$ $ \xch $
\xEH $f(Y) \xcs X \xEd \xCQ $ $ \xch $
\xEP

\xEH
$ \ol{ \ol{T \xcv T' } }= \ol{ \ol{ \ol{T' } } \xcv T}$
\xEH $f(Y \xcs X)=f(Y) \xcs X$
\xEP

\hline

\xEH
$(Log \xFO )$
\xEH $( \xbm \xFO )$
\xEP

\xEH
$ \ol{ \ol{T \xco T' } }$ is one of
\xEH $f(X \xcv Y)$ is one of
\xEP

\xEH
$\ol{\ol{T}},$ or $\ol{\ol{T'}},$ or $\ol{\ol{T}} \xcs \ol{\ol{T'}}$ (by (CCL))
\xEH $f(X),$ $f(Y)$ or $f(X) \xcv f(Y)$
\xEP

\hline

\xEH
$(Log \xcv )$
\xEH $( \xbm \xcv )$
\xEP

\xEH
$Con( \ol{ \ol{T' } } \xcv T),$ $ \xCN Con( \ol{ \ol{T' } }
\xcv \ol{ \ol{T} })$ $ \xch $
\xEH $f(Y) \xcs (X-f(X)) \xEd \xCQ $ $ \xch $
\xEP

\xEH
$ \xCN Con( \ol{ \ol{T \xco T' } } \xcv T' )$
\xEH $f(X \xcv Y) \xcs Y= \xCQ$
\xEP

\hline

\xEH
$(Log \xcv ' )$
\xEH $( \xbm \xcv ' )$
\xEP

\xEH
$Con( \ol{ \ol{T' } } \xcv T),$ $ \xCN Con( \ol{ \ol{T' }
} \xcv \ol{ \ol{T} })$ $ \xch $
\xEH $f(Y) \xcs (X-f(X)) \xEd \xCQ $ $ \xch $
\xEP

\xEH
$ \ol{ \ol{T \xco T' } }= \ol{ \ol{T} }$
\xEH $f(X \xcv Y)=f(X)$
\xEP

\hline

\xEH
\xEH $( \xbm \xbe )$
\xEP

\xEH
\xEH $a \xbe X-f(X)$ $ \xch $
\xEP

\xEH
\xEH $ \xcE b \xbe X.a \xce f(\{a,b\})$
\xEP

\hline

\end{tabular}

}

$(PR)$ is also called infinite conditionalization - we choose the name for
its central role for preferential structures $(PR)$ or $( \xbm PR).$

The system of rules $(AND)$ $(OR)$ $(LLE)$ $(RW)$ $(SC)$ $(CP)$ $(CM)$ $(CUM)$
is also called system $P$ (for preferential), adding $(RatM)$ gives the system
$R$ (for rationality or rankedness).

Roughly: Smooth preferential structures generate logics satisfying system
$P$, ranked structures logics satisfying system $R$.

A logic satisfying $(REF)$, $(ResM)$, and $(CUT)$ is called a consequence
relation.

$(LLE)$ and$(CCL)$ will hold automatically, whenever we work with model sets.

$(AND)$ is obviously closely related to filters, and corresponds to closure
under finite intersections. $(RW)$ corresponds to upward closure of filters.

More precisely, validity of both depend on the definition, and the
direction we consider.

Given $f$ and $(\xbm \xcc )$, $f(X)\xcc X$ generates a pricipal filter:
$\{X'\xcc X:f(X)\xcc X'\}$, with
the definition: If $X=M(T)$, then $T\xcn \xbf$  iff $f(X)\xcc M(\xbf )$.
Validity of $(AND)$ and
$(RW)$ are then trivial.

Conversely, we can define for $X=M(T)$

$\xdx:=\{X'\xcc X: \xcE \xbf (X'=X\xcs M(\xbf )$ and $T\xcn \xbf )\}$.

$(AND)$ then makes $\xdx$  closed under
finite intersections, $(RW)$ makes $\xdx$  upward
closed. This is in the infinite case usually not yet a filter, as not all
subsets of $X$ need to be definable this way.
In this case, we complete $\xdx$  by
adding all $X''$ such that there is $X'\xcc X''\xcc X$, $X'\xbe\xdx$.

Alternatively, we can define

$\xdx:=\{X'\xcc X: \xcS\{X \xcs M(\xbf ): T\xcn \xbf \} \xcc X' \}$.

$(SC)$ corresponds to the choice of a subset.

$(CP)$ is somewhat delicate, as it presupposes that the chosen model set is
non-empty. This might fail in the presence of ever better choices, without
ideal ones; the problem is addressed by the limit versions.

$(PR)$ is an infinitary version of one half of the deduction theorem: Let $T$
stand for $\xbf$, $T'$ for $\xbq$, and $\xbf \xcu \xbq \xcn \xbs$,
so $\xbf \xcn \xbq \xcp \xbs$, but $(\xbq \xcp \xbs )\xcu \xbq \xcl \xbs$.

$(CUM)$ (whose most interesting half in our context is $(CM)$) may best be seen
as
normal use of lemmas: We have worked hard and found some lemmas. Now
we can take a rest, and come back again with our new lemmas. Adding them to the
axioms will neither add new theorems, nor prevent old ones to hold.

\index{Fact Mu-Base}

\ed

\bfa

$\hspace{0.01em}$


\label{Fact Mu-Base}

This table is to be read as follows: If the left hand side holds for some
function $f: \xdy \xcp \xdp (U),$ and the auxiliary properties noted in
the middle also
hold for $f$ or $ \xdy,$ then the right hand side will hold, too - and
conversely.

{\small

\begin{tabular}{|c|c|c|c|}

\hline

\multicolumn{4}{|c|}{Basics} \xEP

\hline

(1.1)
\xEH
$(\xbm PR)$
\xEH
$\xch$ $(\xcs)+(\xbm \xcc)$
\xEH
$(\xbm PR')$
\xEP

\cline{1-1}

\cline{3-3}

(1.2)
\xEH
\xEH
$\xci$
\xEH
\xEP

\hline

(2.1)
\xEH
$(\xbm PR)$
\xEH
$\xch$ $(\xbm \xcc)$
\xEH
$(\xbm OR)$
\xEP

\cline{1-1}

\cline{3-3}

(2.2)
\xEH
\xEH
$\xci$ $(\xbm \xcc)$ + closure
\xEH
\xEP

\xEH
\xEH
under set difference
\xEH
\xEP

\hline

(3)
\xEH
$(\xbm PR)$
\xEH
$\xch$
\xEH
$( \xbm CUT)$
\xEP

\hline

(4)
\xEH
$(\xbm \xcc )+(\xbm \xcc \xcd )+(\xbm CUM)+$
\xEH
$\xcH$
\xEH
$( \xbm PR)$
\xEP

\xEH
$(\xbm RatM)+(\xcs )$
\xEH
\xEH
\xEP

\hline

\multicolumn{4}{|c|}{Cumulativity} \xEP

\hline

(5.1)
\xEH
$(\xbm CM)$
\xEH
$\xch$ $(\xcs)+(\xbm \xcc)$
\xEH
$(\xbm ResM)$
\xEP

\cline{1-1}

\cline{3-3}

(5.2)
\xEH
\xEH
$\xci$ (infin.)
\xEH
\xEP

\hline

(6)
\xEH
$(\xbm CM)+(\xbm CUT)$
\xEH
$\xcj$
\xEH
$(\xbm CUM)$
\xEP

\hline

(7)
\xEH
$( \xbm \xcc )+( \xbm \xcc \xcd )$
\xEH
$\xch$
\xEH
$( \xbm CUM)$
\xEP

\hline

(8)
\xEH
$( \xbm \xcc )+( \xbm CUM)+( \xcs )$
\xEH
$\xch$
\xEH
$( \xbm \xcc \xcd )$
\xEP

\hline

(9)
\xEH
$( \xbm \xcc )+( \xbm CUM)$
\xEH
$\xcH$
\xEH
$( \xbm \xcc \xcd )$
\xEP

\hline

\multicolumn{4}{|c|}{Rationality} \xEP

\hline

(10)
\xEH
$( \xbm RatM )+( \xbm PR )$
\xEH
$\xch$
\xEH
$( \xbm =)$
\xEP

\hline

(11)
\xEH
$( \xbm =)$
\xEH
$ \xch $
\xEH
$( \xbm PR),$
\xEP

\hline

(12.1)
\xEH
$( \xbm =)$
\xEH
$ \xch $ $(\xcs)+( \xbm \xcc )$
\xEH
$( \xbm =' ),$
\xEP
\cline{1-1}
\cline{3-3}
(12.2)
\xEH
\xEH
$ \xci $
\xEH
\xEP

\hline

(13)
\xEH
$( \xbm \xcc ),$ $( \xbm =)$
\xEH
$ \xch $ $(\xcv)$
\xEH
$( \xbm \xcv ),$
\xEP

\hline

(14)
\xEH
$( \xbm \xcc ),$ $( \xbm \xCQ ),$ $( \xbm =)$
\xEH
$ \xch $ $(\xcv)$
\xEH
$( \xbm \xFO ),$ $( \xbm \xcv ' ),$ $( \xbm CUM),$
\xEP

\hline

(15)
\xEH
$( \xbm \xcc )+( \xbm \xFO )$
\xEH
$ \xch $ $\xdy$ closed under set difference
\xEH
$( \xbm =),$
\xEP

\hline

(16)
\xEH
$( \xbm \xFO )+( \xbm \xbe )+( \xbm PR)+$
\xEH
$ \xch $ $(\xcv)$ + $\xdy$ contains singletons
\xEH
$( \xbm =),$
\xEP
\xEH
$( \xbm \xcc )$
\xEH
\xEH
\xEP

\hline

(17)
\xEH
$( \xbm CUM)+( \xbm =)$
\xEH
$ \xch $ $(\xcv)$ + $\xdy$ contains singletons
\xEH
$( \xbm \xbe ),$
\xEP

\hline

(18)
\xEH
$( \xbm CUM)+( \xbm =)+( \xbm \xcc )$
\xEH
$ \xch $ $(\xcv)$
\xEH
$( \xbm \xFO ),$
\xEP

\hline

(19)
\xEH
$( \xbm PR)+( \xbm CUM)+( \xbm \xFO )$
\xEH
$ \xch $ sufficient, e.g. true in $\xdD_{\xdl}$
\xEH
$( \xbm =)$.
\xEP

\hline

(20)
\xEH
$( \xbm \xcc )+( \xbm PR)+( \xbm =)$
\xEH
$ \xcH $
\xEH
$( \xbm \xFO ),$
\xEP

\hline

(21)
\xEH
$( \xbm \xcc )+( \xbm PR)+( \xbm \xFO )$
\xEH
$ \xcH $ (without closure
\xEH
$( \xbm =)$
\xEP
\xEH
\xEH
under set difference),
\xEH
\xEP

\hline

(22)
\xEH
$( \xbm \xcc )+( \xbm PR)+( \xbm \xFO )+$
\xEH
$ \xcH $
\xEH
$( \xbm \xbe )$
\xEP
\xEH
$( \xbm =)+( \xbm \xcv )$
\xEH
\xEH
(thus not representability
\xEP
\xEH
\xEH
\xEH
by ranked structures)
\xEP

\hline

\end{tabular}

}

\index{Proposition Alg-Log}

\efa

\bp

$\hspace{0.01em}$


\label{Proposition Alg-Log}

The following table is to be read as follows:

Let a logic $ \xcn $ satisfies $ \xCf (LLE)$ and $ \xCf (CCL),$ and define
a function $f: \xdD_{ \xdl } \xcp \xdD_{ \xdl }$
by $f(M(T)):=M( \ol{ \ol{T} }).$ Then $f$ is well defined, satisfies $(
\xbm dp),$ and $ \ol{ \ol{T} }=Th(f(M(T))).$

If $ \xcn $ satisfies a rule in the left hand side,
then - provided the additional properties noted in the middle for $ \xch $
hold, too -
$f$ will satisfy the property in the right hand side.

Conversely, if $f: \xdy \xcp \xdp (M_{ \xdl })$ is a function, with $
\xdD_{ \xdl } \xcc \xdy,$ and we define a logic
$ \xcn $ by $ \ol{ \ol{T} }:=Th(f(M(T))),$ then $ \xcn $ satisfies $ \xCf
(LLE)$ and $ \xCf (CCL).$
If $f$ satisfies $( \xbm dp),$ then $f(M(T))=M( \ol{ \ol{T} }).$

If $f$ satisfies a property in the right hand side,
then - provided the additional properties noted in the middle for $ \xci $
hold, too -
$ \xcn $ will satisfy the property in the left hand side.

If ``formula'' is noted in the table, this means that, if one of the
theories
(the one named the same way in Definition \ref{Definition Log-Cond})
is equivalent to a formula, we can renounce on $( \xbm dp).$

{\small

\begin{tabular}{|c|c|c|c|}

\hline

\multicolumn{4}{|c|}{Basics} \xEP

\hline

(1.1) \xEH $(OR)$ \xEH $\xch$ \xEH $(\xbm OR)$ \xEP

\cline{1-1}

\cline{3-3}

(1.2) \xEH \xEH $\xci$ \xEH \xEP

\hline

(2.1) \xEH $(disjOR)$ \xEH $\xch$ \xEH $(\xbm disjOR)$ \xEP

\cline{1-1}

\cline{3-3}

(2.2) \xEH \xEH $\xci$ \xEH \xEP

\hline

(3.1) \xEH $(wOR)$ \xEH $\xch$ \xEH $(\xbm wOR)$ \xEP

\cline{1-1}

\cline{3-3}

(3.2) \xEH \xEH $\xci$ \xEH \xEP

\hline

(4.1) \xEH $(SC)$ \xEH $\xch$ \xEH $(\xbm \xcc)$ \xEP

\cline{1-1}

\cline{3-3}

(4.2) \xEH \xEH $\xci$ \xEH \xEP

\hline

(5.1) \xEH $(CP)$ \xEH $\xch$ \xEH $(\xbm \xCQ)$ \xEP

\cline{1-1}

\cline{3-3}

(5.2) \xEH \xEH $\xci$ \xEH \xEP

\hline

(6.1) \xEH $(PR)$ \xEH $\xch$ \xEH $(\xbm PR)$ \xEP

\cline{1-1}

\cline{3-3}

(6.2) \xEH \xEH $\xci$ $(\xbm dp)+(\xbm \xcc)$ \xEH \xEP

\cline{1-1}

\cline{3-3}

(6.3) \xEH \xEH $\xcI$ without $(\xbm dp)$ \xEH \xEP

\cline{1-1}

\cline{3-3}

(6.4) \xEH \xEH $\xci$ $(\xbm \xcc)$ \xEH \xEP

\xEH \xEH $T'$ a formula \xEH \xEP

\hline

(6.5) \xEH $(PR)$ \xEH $\xci$ \xEH $(\xbm PR')$ \xEP

\xEH \xEH $T'$ a formula \xEH \xEP

\hline

(7.1) \xEH $(CUT)$ \xEH $\xch$ \xEH $(\xbm CUT)$ \xEP

\cline{1-1}

\cline{3-3}

(7.2) \xEH \xEH $\xci$ \xEH \xEP

\hline

\multicolumn{4}{|c|}{Cumulativity} \xEP

\hline

(8.1) \xEH $(CM)$ \xEH $\xch$ \xEH $(\xbm CM)$ \xEP

\cline{1-1}

\cline{3-3}

(8.2) \xEH \xEH $\xci$ \xEH \xEP

\hline

(9.1) \xEH $(ResM)$ \xEH $\xch$ \xEH $(\xbm ResM)$ \xEP

\cline{1-1}

\cline{3-3}

(9.2) \xEH \xEH $\xci$ \xEH \xEP

\hline

(10.1) \xEH $(\xcc \xcd)$ \xEH $\xch$ \xEH $(\xbm \xcc \xcd)$ \xEP

\cline{1-1}

\cline{3-3}

(10.2) \xEH \xEH $\xci$ \xEH \xEP

\hline

(11.1) \xEH $(CUM)$ \xEH $\xch$ \xEH $(\xbm CUM)$ \xEP

\cline{1-1}

\cline{3-3}

(11.2) \xEH \xEH $\xci$ \xEH \xEP

\hline

\multicolumn{4}{|c|}{Rationality} \xEP

\hline

(12.1) \xEH $(RatM)$ \xEH $\xch$ \xEH $(\xbm RatM)$ \xEP

\cline{1-1}

\cline{3-3}

(12.2) \xEH \xEH $\xci$ $(\xbm dp)$ \xEH \xEP

\cline{1-1}

\cline{3-3}

(12.3) \xEH \xEH $\xcI$ without $(\xbm dp)$ \xEH \xEP

\cline{1-1}

\cline{3-3}

(12.4) \xEH \xEH $\xci$ \xEH \xEP

\xEH \xEH $T$ a formula \xEH \xEP

\hline

(13.1) \xEH $(RatM=)$ \xEH $\xch$ \xEH $(\xbm =)$ \xEP

\cline{1-1}

\cline{3-3}

(13.2) \xEH \xEH $\xci$ $(\xbm dp)$ \xEH \xEP

\cline{1-1}

\cline{3-3}

(13.3) \xEH \xEH $\xcI$ without $(\xbm dp)$ \xEH \xEP

\cline{1-1}

\cline{3-3}

(13.4) \xEH \xEH $\xci$ \xEH \xEP

\xEH \xEH $T$ a formula \xEH \xEP

\hline

(14.1) \xEH $(Log = ')$ \xEH $\xch$ \xEH $(\xbm = ')$ \xEP

\cline{1-1}
\cline{3-3}

(14.2) \xEH \xEH $\xci$ $(\xbm dp)$ \xEH \xEP

\cline{1-1}
\cline{3-3}

(14.3) \xEH \xEH $\xcI$ without $(\xbm dp)$ \xEH \xEP

\cline{1-1}
\cline{3-3}

(14.4) \xEH \xEH $\xci$ $T$ a formula \xEH \xEP

\hline

(15.1) \xEH $(Log \xFO )$ \xEH $\xch$ \xEH $(\xbm \xFO )$ \xEP

\cline{1-1}
\cline{3-3}

(15.2) \xEH \xEH $\xci$ \xEH \xEP

\hline

(16.1)
\xEH
$(Log \xcv )$
\xEH
$\xch$ $(\xbm \xcc)+(\xbm =)$
\xEH
$(\xbm \xcv )$
\xEP

\cline{1-1}
\cline{3-3}

(16.2) \xEH \xEH $\xci$ $(\xbm dp)$ \xEH \xEP

\cline{1-1}
\cline{3-3}

(16.3) \xEH \xEH $\xcI$ without $(\xbm dp)$ \xEH \xEP

\hline

(17.1)
\xEH
$(Log \xcv ')$
\xEH
$\xch$ $(\xbm \xcc)+(\xbm =)$
\xEH
$(\xbm \xcv ')$
\xEP

\cline{1-1}
\cline{3-3}

(17.2) \xEH \xEH $\xci$ $(\xbm dp)$ \xEH \xEP

\cline{1-1}
\cline{3-3}

(17.3) \xEH \xEH $\xcI$ without $(\xbm dp)$ \xEH \xEP

\hline

\end{tabular}

}

\index{Definition Pref-Str}

\ep

\bd

$\hspace{0.01em}$


\label{Definition Pref-Str}

Fix $U \xEd \xCQ,$ and consider arbitrary $X.$
Note that this $X$ has not necessarily anything to do with $U,$ or $ \xdu
$ below.
Thus, the functions $ \xbm_{ \xdm }$ below are in principle functions from
$V$ to $V$ - where $V$
is the set theoretical universe we work in.

(A) Preferential models or structures.

(1) The version without copies:

A pair $ \xdm:=<U, \xeb >$ with $U$ an arbitrary set, and $ \xeb $ an
arbitrary binary relation
is called a preferential model or structure.

(2) The version with copies:

A pair $ \xdm:=< \xdu, \xeb >$ with $ \xdu $ an arbitrary set of pairs,
and $ \xeb $ an arbitrary binary
relation is called a preferential model or structure.

If $<x,i> \xbe \xdu,$ then $x$ is intended to be an element of $U,$ and
$i$ the index of the
copy.

We sometimes also need copies of the relation $ \xeb,$ we will then
replace $ \xeb $
by one or several arrows $ \xba $ attacking non-minimal elements, e.g. $x
\xeb y$ will
be written $ \xba:x \xcp y,$ $<x,i> \xeb <y,i>$ will be written $ \xba
:<x,i> \xcp <y,i>,$ and
finally we might have $< \xba,k>:x \xcp y$ and $< \xba,k>:<x,i> \xcp
<y,i>,$ etc.

(B) Minimal elements, the functions $ \xbm_{ \xdm }$

(1) The version without copies:

Let $ \xdm:=<U, \xeb >,$ and define

$ \xbm_{ \xdm }(X)$ $:=$ $\{x \xbe X:$ $x \xbe U$ $ \xcu $ $ \xCN \xcE x'
\xbe X \xcs U.x' \xeb x\}.$

$ \xbm_{ \xdm }(X)$ is called the set of minimal elements of $X$ (in $
\xdm ).$

(2) The version with copies:

Let $ \xdm:=< \xdu, \xeb >$ be as above. Define

$ \xbm_{ \xdm }(X)$ $:=$ $\{x \xbe X:$ $ \xcE <x,i> \xbe \xdu. \xCN \xcE
<x',i' > \xbe \xdu (x' \xbe X$ $ \xcu $ $<x',i' >' \xeb <x,i>)\}.$

Again, by abuse of language, we say that $ \xbm_{ \xdm }(X)$ is the set of
minimal elements
of $X$ in the structure. If the context is clear, we will also write just
$ \xbm.$

We sometimes say that $<x,i>$ ``kills'' or ``minimizes'' $<y,j>$ if
$<x,i> \xeb <y,j>.$ By abuse of language we also say a set $X$ kills or
minimizes a set
$Y$ if for all $<y,j> \xbe \xdu,$ $y \xbe Y$ there is $<x,i> \xbe \xdu,$
$x \xbe X$ s.t. $<x,i> \xeb <y,j>.$

$ \xdm $ is also called injective or 1-copy, iff there is always at most
one copy
$<x,i>$ for each $x.$ Note that the existence of copies corresponds to a
non-injective labelling function - as is often used in nonclassical
logic, e.g. modal logic.

We say that $ \xdm $ is transitive, irreflexive, etc., iff $ \xeb $ is.

Note that $ \xbm (X)$ might well be empty, even if $X$ is not.
\index{Definition Pref-Log}

\ed

\bd

$\hspace{0.01em}$


\label{Definition Pref-Log}

We define the consequence relation of a preferential structure for a
given propositional language $ \xdl.$

(A)

(1) If $m$ is a classical model of a language $ \xdl,$ we say by abuse
of language

$<m,i> \xcm \xbf $ iff $m \xcm \xbf,$

and if $X$ is a set of such pairs, that

$X \xcm \xbf $ iff for all $<m,i> \xbe X$ $m \xcm \xbf.$

(2) If $ \xdm $ is a preferential structure, and $X$ is a set of $ \xdl
-$models for a
classical propositional language $ \xdl,$ or a set of pairs $<m,i>,$
where the $m$ are
such models, we call $ \xdm $ a classical preferential structure or model.

(B)

Validity in a preferential structure, or the semantical consequence
relation
defined by such a structure:

Let $ \xdm $ be as above.

We define:

$T \xcm_{ \xdm } \xbf $ iff $ \xbm_{ \xdm }(M(T)) \xcm \xbf,$ i.e. $
\xbm_{ \xdm }(M(T)) \xcc M( \xbf ).$

$ \xdm $ will be called definability preserving iff for all $X \xbe \xdD_{
\xdl }$ $ \xbm_{ \xdm }(X) \xbe \xdD_{ \xdl }.$

As $ \xbm_{ \xdm }$ is defined on $ \xdD_{ \xdl },$ but need by no means
always result in some new
definable set, this is (and reveals itself as a quite strong) additional
property.
\index{Definition Smooth}

\ed

\bd

$\hspace{0.01em}$


\label{Definition Smooth}

Let $ \xdy \xcc \xdp (U).$ (In applications to logic, $ \xdy $ will be $
\xdD_{ \xdl }.)$

A preferential structure $ \xdm $ is called $ \xdy -$smooth iff in every
$X \xbe \xdy $ every element
$x \xbe X$ is either minimal in $X$ or above an element, which is minimal
in $X.$ More
precisely:

(1) The version without copies:

If $x \xbe X \xbe \xdy,$ then either $x \xbe \xbm (X)$ or there is $x'
\xbe \xbm (X).x' \xeb x.$

(2) The version with copies:

If $x \xbe X \xbe \xdy,$ and $<x,i> \xbe \xdu,$ then either there is no
$<x',i' > \xbe \xdu,$ $x' \xbe X,$
$<x',i' > \xeb <x,i>$ or there is $<x',i' > \xbe \xdu,$ $<x',i' > \xeb
<x,i>,$ $x' \xbe X,$ s.t. there is
no $<x'',i'' > \xbe \xdu,$ $x'' \xbe X,$ with $<x'',i'' > \xeb <x',i'
>.$

When considering the models of a language $ \xdl,$ $ \xdm $ will be
called smooth iff
it is $ \xdD_{ \xdl }-$smooth; $ \xdD_{ \xdl }$ is the default.

Obviously, the richer the set $ \xdy $ is, the stronger the condition $
\xdy -$smoothness
will be.
\index{Table Pref-Representation-Without-Ref}

\ed

The following table summarizes representation by not necessarily ranked
preferential structures. The implications on the right are shown in
Proposition \ref{Proposition Alg-Log} (going via the $ \xbm -$functions),
those on the left
are shown in the respective representation theorems.
\label{Table Pref-Representation-Without-Ref}

{\scriptsize

\begin{tabular}{|c|c|c|c|c|}

\hline

$\xbm-$ function
\xEH
\xEH
Pref.Structure
\xEH
\xEH
Logic
\xEP

\hline

$(\xbm \xcc)+(\xbm PR)$
\xEH
$\xci$
\xEH
general
\xEH
$\xch$ $(\xbm dp)$
\xEH
$(LLE)+(RW)+(SC)+(PR)$
\xEP

\cline{2-2}
\cline{4-4}

\xEH
$\xch$
\xEH
\xEH
$\xci$
\xEH
\xEP

\cline{2-2}
\cline{4-4}

\xEH
\xEH
\xEH
$\xcH$ without $(\xbm dp)$
\xEH
\xEP

\hline

$(\xbm \xcc)+(\xbm PR)$
\xEH
$\xci$
\xEH
transitive
\xEH
$\xch$ $(\xbm dp)$
\xEH
$(LLE)+(RW)+(SC)+(PR)$
\xEP

\cline{2-2}
\cline{4-4}

\xEH
$\xch$
\xEH
\xEH
$\xci$
\xEH
\xEP

\cline{2-2}
\cline{4-4}

\xEH
\xEH
\xEH
$\xcH$ without $(\xbm dp)$
\xEH
\xEP

\hline

$(\xbm \xcc)+(\xbm PR)+(\xbm CUM)$
\xEH
$\xci$
\xEH
smooth
\xEH
$\xch$ $(\xbm dp)$
\xEH
$(LLE)+(RW)+(SC)+(PR)+$
\xEP

\xEH
\xEH
\xEH
\xEH
$(CUM)$
\xEP

\cline{2-2}
\cline{4-4}

\xEH
$\xch$ $(\xcv)$
\xEH
\xEH
$\xci$ $(\xcv)$
\xEH
\xEP

\cline{2-2}
\cline{4-4}

\xEH
\xEH
\xEH
$\xcH$ without $(\xbm dp)$
\xEH
\xEP

\hline

$(\xbm \xcc)+(\xbm PR)+(\xbm CUM)$
\xEH
$\xci$
\xEH
smooth+transitive
\xEH
$\xch$ $(\xbm dp)$
\xEH
$(LLE)+(RW)+(SC)+(PR)+$
\xEP

\xEH
\xEH
\xEH
\xEH
$(CUM)$
\xEP

\cline{2-2}
\cline{4-4}

\xEH
$\xch$ $(\xcv)$
\xEH
\xEH
$\xci$ $(\xcv)$
\xEH
\xEP

\cline{2-2}
\cline{4-4}

\xEH
\xEH
\xEH
$\xcH$ without $(\xbm dp)$
\xEH
\xEP

\hline

\end{tabular}

}

\section{Motivation - two sequent calculi}
\label{Section Sequent-Calculi}
\subsection{Introduction}
\label{Section Sequent Calculi Introduction}

This section serves mainly as a posteriori motivation for our examination
of weak closure conditions of the domain. The second author realized first
when
looking at Lehmann's plausibility logic, that absence of $( \xcv )$ might
be a
problem for representation.

Beyond motivation, the reader will see here two ``real life'' examples where
closure under $( \xcv )$ is not given, and thus problems arise. So this is
also
a warning against a too naive treatment of representation problems,
neglecting
domain closure issues.
\index{Plausibility Logic}
\subsection{Plausibility Logic}
\label{Section Plausibility Logic}

\paragraph{
Discussion of plausibility logic
}

$\hspace{0.01em}$


\label{Section Discussion of plausibility logic}

Plausibility logic was introduced by $D.$ Lehmann  \cite{Leh92a},
 \cite{Leh92b}
as a sequent
calculus in a propositional language without connectives. Thus, a
plausibility
logic language $ \xdl $ is just a set, whose elements correspond to
propositional
variables, and a sequent has the form $X \xcn Y,$ where $X,$ $Y$ are $
\ul{finite}$ subsets of
$ \xdl,$ thus, in the intuitive reading, $ \xcU X \xcn \xcO Y.$ (We use $
\xcn $ instead of the $ \xcl $ used
in  \cite{Leh92a},  \cite{Leh92b} and continue to reserve $
\xcl $ for classical logic.)

\paragraph{
The details:
}

$\hspace{0.01em}$


\label{Section The details:}

\bn

$\hspace{0.01em}$


\label{Notation Plausi-1}

We abuse notation, and write $X \xcn a$ for $X \xcn \{a\},$ $X,a \xcn Y$
for $X \xcv \{a\} \xcn Y,$ $ab \xcn Y$ for
$\{a,b\} \xcn Y,$ etc. When discussing plausibility logic, $X,Y,$ etc.
will denote finite
subsets of $ \xdl,$ $a,b,$ etc. elements of $ \xdl.$

We first define the logical properties we will examine.

\en

\bd

$\hspace{0.01em}$


\label{Definition Plausi-1}

$X$ and $Y$ will be finite subsets of $ \xdl,$ a, etc. elements of $ \xdl
.$
The base axiom and rules of plausibility logic are
(we use the prefix ``Pl'' to differentiate them from the usual ones):

(PlI) (Inclusion): $X \xcn a$ for all $a \xbe X,$

(PlRM) (Right Monotony): $X \xcn Y$ $ \xch $ $X \xcn a,Y,$

(PlCLM) (Cautious Left Monotony): $X \xcn a,$ $X \xcn Y$ $ \xch $ $X,a
\xcn Y,$

(PlCC) (Cautious Cut): $X,a_{1} \Xl a_{n} \xcn Y,$ and for all $1 \xck i
\xck n$ $X \xcn a_{i},Y$ $ \xch $ $X \xcn Y,$

and as a special case of (PlCC):

(PlUCC) (Unit Cautious Cut): $X,a \xcn Y$, $X \xcn a,Y$ $ \xch $ $X \xcn
Y.$

and we denote by PL, for plausibility logic, the full system, i.e.
$(PlI)+(PlRM)+(PlCLM)+(PlCC).$ $ \xcz $
\\[3ex]

\ed

We now adapt the definition of a preferential model to plausibility logic.
This is the central definition on the semantic side.

\bd

$\hspace{0.01em}$


\label{Definition Plausi-2}

Fix a plausibility logic language $ \xdl.$ A model for $ \xdl $ is then
just an arbitrary
subset of $ \xdl.$

If $ \xdm:=<M, \xeb >$ is a preferential model s.t. $M$ is a set of
(indexed) $ \xdl -$models,
then for a finite set $X \xcc \xdl $ (to be imagined on the left hand side
of $ \xcn $!), we
define

(a) $m \xcm X$ iff $X \xcc m$

(b) $M(X)$ $:=$ $\{m$: $<m,i> \xbe M$ for some $i$ and $m \xcm X\}$

(c) $ \xbm (X)$ $:=$ $\{m \xbe M(X)$: $ \xcE <m,i> \xbe M. \xCN \xcE <m'
,i' > \xbe M$ $(m' \xbe M(X)$ $ \xcu $ $<m',i' > \xeb <m,i>)\}$

(d) $X \xcm_{ \xdm }Y$ iff $ \xcA m \xbe \xbm (X).m \xcs Y \xEd \xCQ.$ $
\xcz $
\\[3ex]

\ed

(a) reflects the intuitive reading of $X$ as $ \xcU X,$ and (d) that of
$Y$ as $ \xcO Y$ in
$X \xcn Y.$ Note that $X$ is a set of ``formulas'', and $ \xbm (X)= \xbm_{
\xdm }(M(X)).$

We note as trivial consequences of the definition.

\bfa

$\hspace{0.01em}$


\label{Fact Plausi-1}

(a) $a \xcm_{ \xdm }b$ iff for all $m \xbe \xbm (a).b \xbe m$

(b) $X \xcm_{ \xdm }Y$ iff $ \xbm (X) \xcc \xcV \{M(b):b \xbe Y\}$

(c) $m \xbe \xbm (X)$ $ \xcu $ $X \xcc X' $ $ \xcu $ $m \xbe M(X' )$ $
\xcp $ $m \xbe \xbm (X' )$. $ \xcz $
\\[3ex]

\efa

We note without proof: $(PlI)+(PlRM)+(PlCC)$ is complete (and sound) for
preferential models

We note the following fact for smooth preferential models:

\bfa

$\hspace{0.01em}$


\label{Fact Plausi-2}

Let $ \xCf U,X,Y$ be any sets, $ \xdm $ be smooth for at least $\{Y,X\}$
and
let $ \xbm (Y) \xcc U \xcv X,$ $ \xbm (X) \xcc U,$ then $X \xcs Y \xcs
\xbm (U) \xcc \xbm (Y).$ (This is, of course,
a special case of $( \xbm Cum1),$ see Definition \ref{Definition Cum-Alpha}.

\efa

\be

$\hspace{0.01em}$


\label{Example Plausi-1}

Let $ \xdl:=\{a,b,c,d,e,f\},$ and
$ \xdx $ $:=$ $\{a \xcn b$, $b \xcn a$, $a \xcn c$, $a \xcn fd$, $dc
\xcn ba$, $dc \xcn e$, $fcba \xcn e\}.$
We show that $ \xdx $ does not have a smooth representation.

\ee

\bfa

$\hspace{0.01em}$


\label{Fact Plausi-3}

$ \xdx $ does not entail $a \xcn e.$

\efa

See  \cite{Sch96-3} for a proof.

Suppose now that there is a smooth preferential model $ \xdm =<M, \xeb >$
for plausibility
logic which represents $ \xcn,$ i.e. for all $ \xCf X,Y$ finite subsets
of $ \xdl $ $X \xcn Y$ iff
$X \xcm_{ \xdm }Y.$ (See Definition \ref{Definition Plausi-2} and Fact \ref{Fact
Plausi-1}.)

$a \xcn a,$ $a \xcn b,$ $a \xcn c$ implies for $m \xbe \xbm (a)$ $a,b,c
\xbe m.$ Moreover, as $a \xcn df,$ then also
$d \xbe m$ or $f \xbe m.$ As $a \xcN e,$ there must be $m \xbe \xbm (a)$
s.t. $e \xce m.$ Suppose now $m \xbe \xbm (a)$
with $f \xbe m.$ So $a,b,c,f \xbe m,$ thus by $m \xbe \xbm (a)$ and Fact
\ref{Fact Plausi-1},
$m \xbe \xbm (a,b,c,f).$ But
$fcba \xcn e,$ so $e \xbe m.$ We thus have shown that $m \xbe \xbm (a)$
and $f \xbe m$ implies $e \xbe m.$
Consequently, there must be $m \xbe \xbm (a)$ s.t. $d \xbe m,$ $e \xce m.$
Thus, in particular, as $cd \xcn e,$ there is $m \xbe \xbm (a),$ $a,b,c,d
\xbe m,$ $m \xce \xbm (cd).$
But by $cd \xcn ab,$ and $b \xcn a,$ $ \xbm (cd) \xcc M(a) \xcv M(b)$ and
$ \xbm (b) \xcc M(a)$ by
Fact \ref{Fact Plausi-1}.
Let now $T:=M(cd),$ $R:=M(a),$ $S:=M(b),$ and $ \xbm_{ \xdm }$ be the
choice function of the
minimal elements in the structure $ \xdm,$ we then have by $ \xbm (S)=
\xbm_{ \xdm }(M(S))$:

1. $ \xbm_{ \xdm }(T) \xcc R \xcv S,$

2. $ \xbm_{ \xdm }(S) \xcc R,$

3. there is $m \xbe S \xcs T \xcs \xbm_{ \xdm }(R),$ but $m \xce \xbm_{
\xdm }(T),$

but this contradicts above Fact \ref{Fact Plausi-2}. $ \xcz $
(Example Plausi-1)
\\[3ex]
\index{Arieli-Avron}
\subsection{A comment on the work by Arieli and Avron}
\label{Section Arieli-Avron}

We turn to a similar case, published in  \cite{AA00}.
Definitions are due to  \cite{AA00}, for motivation the reader is
referred there.

\bd

$\hspace{0.01em}$


\label{Definition Arieli-Avron-1}

(1) A Scott consequence relation, abbreviated scr, is a binary relation $
\xcl $
between sets of formulae, that satisfies the following conditions:

(s-R) if $ \xbG \xcs \xbD \xEd \xCQ,$ the $ \xbG \xcl \xbD $
(M) if $ \xbG \xcl \xbD $ and $ \xbG \xcc \xbG ',$ $ \xbD \xcc \xbD ',$
then $ \xbG ' \xcl \xbD ' $
(C) if $ \xbG \xcl \xbq, \xbD $ and $ \xbG ', \xbq \xcl \xbD ',$ then $
\xbG, \xbG ' \xcl \xbD, \xbD ' $

(2) A Scott cautious consequence relation, abbreviated sccr, is a binary
relation $ \xcn $ between nonempty sets of formulae, that satisfies the
following
conditions:

(s-R) if $ \xbG \xcs \xbD \xEd \xCQ,$ the $ \xbG \xcn \xbD $
(CM) if $ \xbG \xcn \xbD $ and $ \xbG \xcn \xbq,$ then $ \xbG, \xbq \xcn
\xbD $
(CC) if $ \xbG \xcn \xbq $ and $ \xbG, \xbq \xcn \xbD,$ then $ \xbG \xcn
\xbD.$

\ed

\be

$\hspace{0.01em}$


\label{Example Arieli-Avron-1}

We have two consequence relations, $ \xcl $ and $ \xcn.$

The rules to consider are

$LCC^{n}$ $ \frac{ \xbG \xcn \xbq_{1}, \xbD  \Xl  \xbG \xcn \xbq_{n}, \xbD
\xbG, \xbq_{1}, \Xl, \xbq_{n} \xcn \xbD }{ \xbG \xcn \xbD }$

$RW^{n}$ $ \frac{ \xbG \xcn \xbq_{i}, \xbD i=1 \Xl n \xbG, \xbq_{1}, \Xl
, \xbq_{n} \xcl \xbf }{ \xbG \xcn \xbf, \xbD }$

Cum $ \xbG, \xbD \xEd \xCQ,$ $ \xbG \xcl \xbD $ $ \xcp $ $ \xbG \xcn
\xbD $

RM $ \xbG \xcn \xbD $ $ \xcp $ $ \xbG \xcn \xbq, \xbD $

CM $ \frac{ \xbG \xcn \xbq \xbG \xcn \xbD }{ \xbG, \xbq \xcn \xbD }$

$s-R$ $ \xbG \xcs \xbD \xEd \xCQ $ $ \xcp $ $ \xbG \xcn \xbD $

$M$ $ \xbG \xcl \xbD,$ $ \xbG \xcc \xbG ',$ $ \xbD \xcc \xbD ' $ $ \xcp
$ $ \xbG ' \xcl \xbD ' $

$C$ $ \frac{ \xbG_{1} \xcl \xbq, \xbD_{1} \xbG_{2}, \xbq \xcl \xbD_{2}}{
\xbG_{1}, \xbG_{2} \xcl \xbD_{1}, \xbD_{2}}$

Let $ \xdl $ be any set.
Define now $ \xbG \xcl \xbD $ iff $ \xbG \xcs \xbD \xEd \xCQ.$
Then $s-R$ and $M$ for $ \xcl $ are trivial. For $C:$ If $ \xbG_{1} \xcs
\xbD_{1} \xEd \xCQ $ or $ \xbG_{1} \xcs \xbD_{1} \xEd \xCQ,$ the
result is trivial. If not, $ \xbq \xbe \xbG_{1}$ and $ \xbq \xbe
\xbD_{2},$ which implies the result.
So $ \xcl $ is a scr.

Consider now the rules for a sccr which is $ \xcl -$plausible for this $
\xcl.$
Cum is equivalent to $s-$R, which is essentially (PlI) of Plausibility
Logic.
Consider $RW^{n}.$ If $ \xbf $ is one of the $ \xbq_{i},$ then the
consequence $ \xbG \xcn \xbf, \xbD $ is a case
of one of the other hypotheses. If not, $ \xbf \xbe \xbG,$ so $ \xbG \xcn
\xbf $ by $s-$R, so $ \xbG \xcn \xbf, \xbD $
by RM (if $ \xbD $ is finite). So, for this $ \xcl,$ $RW^{n}$ is a
consequence of $s-R$ $+$ RM.

We are left with $LCC^{n},$ RM, CM, $s-$R, it was shown in  \cite{Sch04}
and  \cite{Sch96-3}
that this does not suffice to guarantee smooth representability, by
failure of
$( \xbm Cum1)$ - see Definition \ref{Definition Cum-Alpha}.
\section{Cumulativity without $(\xcv)$}
\label{Section Cumulativity-without-union}
\subsection{Introduction}
\label{Section Cumulativity-without-union Introduction}
\index{Comment Cum-Union}

\ee

\bcom

$\hspace{0.01em}$


\label{Comment Cum-Union}

We show here that, without sufficient closure
properties, there is an infinity of versions of cumulativity, which
collaps
to usual cumulativity when the domain is closed under finite unions.
Closure properties thus reveal themselves as a powerful tool to show
independence of properties.

We work in some fixed arbitrary set $Z,$ all sets considered will be
subsets of $Z.$

Unless said otherwise, we use without further
mentioning $( \xbm PR)$ and $( \xbm \xcc ).$

\ecom

We first give the definition of the new conditions in this introduction,
and then state and prove the main
result (Example \ref{Example Inf-Cum-Alpha}), and
finally show some properties of the new
conditions (Fact \ref{Fact Cum-Alpha}).
\index{Definition Cum-Alpha}

\bd

$\hspace{0.01em}$


\label{Definition Cum-Alpha}

For any ordinal $ \xba,$ we define

$( \xbm Cum \xba ):$

If for all $ \xbb \xck \xba $ $ \xbm (X_{ \xbb }) \xcc U \xcv \xcV \{X_{
\xbg }: \xbg < \xbb \}$ hold, then so does
$ \xcS \{X_{ \xbg }: \xbg \xck \xba \} \xcs \xbm (U) \xcc \xbm (X_{ \xba
}).$

$( \xbm Cumt \xba ):$

If for all $ \xbb \xck \xba $ $ \xbm (X_{ \xbb }) \xcc U \xcv \xcV \{X_{
\xbg }: \xbg < \xbb \}$ hold, then so does
$X_{ \xba } \xcs \xbm (U) \xcc \xbm (X_{ \xba }).$

( `` $t$ '' stands for transitive, see Fact \ref{Fact Cum-Alpha}, (2.2)
below.)

$( \xbm Cum \xca )$ and $( \xbm Cumt \xca )$ will be the class of all $(
\xbm Cum \xba )$ or $( \xbm Cumt \xba )$ -
read their ``conjunction'', i.e. if we say that $( \xbm Cum \xca )$ holds,
we mean that
all $( \xbm Cum \xba )$ hold.
\index{Note Cum-Alpha}

\ed

\paragraph{
Note Cum-Alpha
}

$\hspace{0.01em}$


\label{Section Note Cum-Alpha}

The first conditions thus have the form:

$( \xbm Cum0)$ $ \xbm (X_{0}) \xcc U$ $ \xcp $ $X_{0} \xcs \xbm (U) \xcc
\xbm (X_{0}),$

$( \xbm Cum1)$ $ \xbm (X_{0}) \xcc U,$ $ \xbm (X_{1}) \xcc U \xcv X_{0}$ $
\xcp $ $X_{0} \xcs X_{1} \xcs \xbm (U) \xcc \xbm (X_{1}),$

$( \xbm Cum2)$ $ \xbm (X_{0}) \xcc U,$ $ \xbm (X_{1}) \xcc U \xcv X_{0},$
$ \xbm (X_{2}) \xcc U \xcv X_{0} \xcv X_{1}$ $ \xcp $ $X_{0} \xcs X_{1}
\xcs X_{2} \xcs \xbm (U) \xcc \xbm (X_{2}).$

$( \xbm Cumt \xba )$ differs from $( \xbm Cum \xba )$ only in the
consequence, the intersection contains
only the last $X_{ \xba }$ - in particular, $( \xbm Cum0)$ and $( \xbm
Cumt0)$ coincide.

Recall that condition $( \xbm Cum1)$ is the crucial condition in  \cite{Leh92a},
which
failed, despite $( \xbm CUM),$ but which has to hold in all smooth models.
This
condition $( \xbm Cum1)$ was the starting point of the investigation.
\subsection{The results}
\label{Section Cumulativity-without-union Results}
\index{Example Inf-Cum-Alpha}

\be

$\hspace{0.01em}$


\label{Example Inf-Cum-Alpha}

This important example shows that the conditions $( \xbm Cum \xba )$ and
$( \xbm Cumt \xba )$ defined in Definition \ref{Definition Cum-Alpha}
are all different in the absence of $( \xcv ),$ in its
presence they all collaps (see Fact \ref{Fact Cum-Alpha} below). More
precisely,
the following (class of) $example(s)$ shows that the $( \xbm Cum \xba )$
increase in
strength. For any finite or infinite ordinal $ \xbk >0$ we construct an
example s.t.

(a) $( \xbm PR)$ and $( \xbm \xcc )$ hold

(b) $( \xbm CUM)$ holds

(c) $( \xcS )$ holds

(d) $( \xbm Cumt \xba )$ holds for $ \xba < \xbk $

(e) $( \xbm Cum \xbk )$ fails.

\ee

\subparagraph{
Proof:
}

$\hspace{0.01em}$


\label{Section Proof:}

We define a suitable base set and a non-transitive binary relation $ \xeb
$
on this set, as well as a suitable set $ \xdx $ of subsets, closed under
arbitrary
intersections, but not under finite unions, and define $ \xbm $ on these
subsets
as usual in preferential structures by $ \xeb.$ Thus, $( \xbm PR)$ and $(
\xbm \xcc )$ will hold.
It will be immediate that $( \xbm Cum \xbk )$ fails, and we will show that
$( \xbm CUM)$ and
$( \xbm Cumt \xba )$ for $ \xba < \xbk $ hold by examining the cases.

For simplicity, we first define a set of generators for $ \xdx,$ and
close under
$( \xcS )$ afterwards. The set $U$ will have a special position, it is the
``useful''
starting point to construct chains corresponding to above definitions of
$( \xbm Cum \xba )$ and $( \xbm Cumt \xba ).$

In the sequel,
$i,j$ will be successor ordinals, $ \xbl $ etc. limit ordinals, $ \xba,$
$ \xbb,$ $ \xbk $ any ordinals,
thus e.g. $ \xbl \xck \xbk $ will imply that $ \xbl $ is a limit ordinal $
\xck \xbk,$ etc.

\paragraph{
The base set and the relation $ \xeb $:
}

$\hspace{0.01em}$


\label{Section The base set and the relation  b:}

$ \xbk >0$ is fixed, but arbitrary. We go up to $ \xbk >0.$

The base set is $\{a,b,c\}$ $ \xcv $ $\{d_{ \xbl }: \xbl \xck \xbk \}$ $
\xcv $ $\{x_{ \xba }: \xba \xck \xbk +1\}$ $ \xcv $ $\{x'_{ \xba }: \xba
\xck \xbk \}.$
$a \xeb b \xeb c,$ $x_{ \xba } \xeb x_{ \xba +1},$ $x_{ \xba } \xeb x'_{
\xba },$ $x'_{0} \xeb x_{ \xbl }$ (for any $ \xbl )$ - $ \xeb $ is NOT
transitive.

\paragraph{
The generators:
}

$\hspace{0.01em}$


\label{Section The generators:}

$U:=\{a,c,x_{0}\} \xcv \{d_{ \xbl }: \xbl \xck \xbk \}$ - i.e. $ \Xl
.\{d_{ \xbl }:lim( \xbl ) \xcu \xbl \xck \xbk \},$

$X_{i}:=\{c,x_{i},x'_{i},x_{i+1}\}$ $(i< \xbk ),$

$X_{ \xbl }:=\{c,d_{ \xbl },x_{ \xbl },x'_{ \xbl },x_{ \xbl +1}\} \xcv
\{x'_{ \xba }: \xba < \xbl \}$ $( \xbl < \xbk ),$

$X'_{ \xbk }:=\{a,b,c,x_{ \xbk },x'_{ \xbk },x_{ \xbk +1}\}$ if $ \xbk $
is a successor,

$X'_{ \xbk }:=\{a,b,c,d_{ \xbk },x_{ \xbk },x'_{ \xbk },x_{ \xbk +1}\}
\xcv \{x'_{ \xba }: \xba < \xbk \}$ if $ \xbk $ is a limit.

Thus, $X'_{ \xbk }=X_{ \xbk } \xcv \{a,b\}$ if $X_{ \xbk }$ were defined.

Note that there is only one $X'_{ \xbk },$ and $X_{ \xba }$ is defined
only for $ \xba < \xbk,$ so we will
not have $X_{ \xba }$ and $X'_{ \xba }$ at the same time.

Thus, the values of the generators under $ \xbm $ are:

$ \xbm (U)=U,$

$ \xbm (X_{i})=\{c,x_{i}\},$

$ \xbm (X_{ \xbl })=\{c,d_{ \xbl }\} \xcv \{x'_{ \xba }: \xba < \xbl \},$

$ \xbm (X'_{i})=\{a,x_{i}\}$ $(i>0!),$

$ \xbm (X'_{ \xbl })=\{a,d_{ \xbl }\} \xcv \{x'_{ \xba }: \xba < \xbl \}.$

(We do not assume that the domain is closed under $ \xbm.)$

\paragraph{
Intersections:
}

$\hspace{0.01em}$


\label{Section Intersections:}

We consider first pairwise intersections:

(1) $U \xcs X_{0}=\{c,x_{0}\},$

(2) $U \xcs X_{i}=\{c\},$ $i>0,$

(3) $U \xcs X_{ \xbl }=\{c,d_{ \xbl }\},$

(4) $U \xcs X'_{i}=\{a,c\}$ $(i>0!),$

(5) $U \xcs X'_{ \xbl }=\{a,c,d_{ \xbl }\},$

(6) $X_{i} \xcs X_{j}:$

(6.1) $j=i+1$ $\{c,x_{i+1}\},$

(6.2) else $\{c\},$

(7) $X_{i} \xcs X_{ \xbl }:$

(7.1) $i< \xbl $ $\{c,x'_{i}\},$

(7.2) $i= \xbl +1$ $\{c,x_{ \xbl +1}\},$

(7.3) $i> \xbl +1$ $\{c\},$

(8) $X_{ \xbl } \xcs X_{ \xbl ' }:$ $\{c\} \xcv \{x'_{ \xba }: \xba \xck
min( \xbl, \xbl ' )\}.$

As $X'_{ \xbk }$ occurs only once, $X_{ \xba } \xcs X'_{ \xbk }$ etc. give
no new results.

Note that $ \xbm $ is constant on all these pairwise intersections.

Iterated intersections:

As $c$ is an element of all sets, sets of the type $\{c,z\}$ do not give
any
new results. The possible subsets of $\{a,c,d_{ \xbl }\}:$ $\{c\},$
$\{a,c\},$ $\{c,d_{ \xbl }\}$ exist
already. Thus, the only source of new sets via iterated intersections is
$X_{ \xbl } \xcs X_{ \xbl ' }=\{c\} \xcv \{x'_{ \xba }: \xba \xck min(
\xbl, \xbl ' )\}.$ But, to intersect them, or with some
old sets, will not generate any new sets either. Consequently, the example
satisfies $( \xcS )$ for $ \xdx $ defined by $U,$ $X_{i}$ $(i< \xbk ),$
$X_{ \xbl }$ $( \xbl < \xbk ),$ $X'_{ \xbk },$ and above
paiwise intersections.

We will now verify the positive properties. This is tedious, but
straightforward, we have to check the different cases.

\paragraph{
Validity of $( \xbm CUM)$:
}

$\hspace{0.01em}$


\label{Section Validity of ( mCUM):}

Consider the prerequisite $ \xbm (X) \xcc Y \xcc X.$ If $ \xbm (X)=X$ or
if $X- \xbm (X)$ is a singleton,
$X$ cannot give a violation of $( \xbm CUM).$ So we are left with the
following
candidates for $X:$

(1) $X_{i}:=\{c,x_{i},x'_{i},x_{i+1}\},$ $ \xbm (X_{i})=\{c,x_{i}\}$

Interesting candidates for $Y$ will have 3 elements, but they will all
contain
a. (If $ \xbk < \xbo:$ $U=\{a,c,x_{0}\}.)$

(2) $X_{ \xbl }:=\{c,d_{ \xbl },x_{ \xbl },x'_{ \xbl },x_{ \xbl +1}\} \xcv
\{x'_{ \xba }: \xba < \xbl \},$ $ \xbm (X_{ \xbl })=\{c,d_{ \xbl }\} \xcv
\{x'_{ \xba }: \xba < \xbl \}$

The only sets to contain $d_{ \xbl }$ are $X_{ \xbl },$ $U,$ $U \xcs X_{
\xbl }.$ But $a \xbe U,$ and
$U \xcs X_{ \xbl }$ ist finite. $(X_{ \xbl }$ and $X'_{ \xbl }$ cannot be
present at the same time.)

(3) $X'_{i}:=\{a,b,c,x_{i},x'_{i},x_{i+1}\},$ $ \xbm (X'_{i})=\{a,x_{i}\}$

a is only in $U,$ $X'_{i},$ $U \xcs X'_{i}=\{a,c\},$ but $x_{i} \xce U,$
as $i>0.$

(4) $X'_{ \xbl }:=\{a,b,c,d_{ \xbl },x_{ \xbl },x'_{ \xbl },x_{ \xbl +1}\}
\xcv \{x'_{ \xba }: \xba < \xbl \},$ $ \xbm (X'_{ \xbl })=\{a,d_{ \xbl }\}
\xcv \{x'_{ \xba }: \xba < \xbl \}$

$d_{ \xbl }$ is only in $X'_{ \xbl }$ and $U,$ but $U$ contains no $x'_{
\xba }.$

Thus, $( \xbm CUM)$ holds trivially.

\paragraph{
$( \xbm Cumt \xba )$ hold for $ \xba < \xbk $:
}

$\hspace{0.01em}$


\label{Section ( mCumt a) hold for  a< k:}

To simplify language, we say that we reach $Y$ from $X$ iff $X \xEd Y$ and
there is a
sequence $X_{ \xbb },$ $ \xbb \xck \xba $ and $ \xbm (X_{ \xbb }) \xcc X
\xcv \xcV \{X_{ \xbg }: \xbg < \xbb \},$ and $X_{ \xba }=Y,$ $X_{0}=X.$
Failure of $( \xbm Cumt \xba )$ would then mean that there are $X$ and
$Y,$ we can reach
$Y$ from $X,$ and $x \xbe ( \xbm (X) \xcs Y)- \xbm (Y).$ Thus, in a
counterexample, $Y= \xbm (Y)$ is
impossible, so none of the intersections can be such $Y.$

To reach $Y$ from $X,$ we have to get started from $X,$ i.e. there must be
$Z$ s.t.
$ \xbm (Z) \xcc X,$ $Z \xcC X$ (so $ \xbm (Z) \xEd Z).$ Inspection of the
different cases shows that
we cannot reach any set $Y$ from any case of the intersections, except
from
(1), (6.1), (7.2).

If $Y$ contains a globally minimal element (i.e. there is no smaller
element in
any set), it can only be reached from any $X$ which already contains this
element. The globally minimal elements are a, $x_{0},$ and the $d_{ \xbl
},$ $ \xbl \xck \xbk.$

By these observations, we see that $X_{ \xbl }$ and $X'_{ \xbk }$ can only
be reached from $U.$
From no $X_{ \xba }$ $U$ can be reached, as the globally minimal a is
missing. But
$U$ cannot be reached from $X'_{ \xbk }$ either, as the globally minimal
$x_{0}$ is missing.

When we look at the relation $ \xeb $ defining $ \xbm,$ we see that we
can reach $Y$ from $X$
only by going upwards, adding bigger elements. Thus, from $X_{ \xba },$ we
cannot reach
any $X_{ \xbb },$ $ \xbb < \xba,$ the same holds for $X'_{ \xbk }$ and
$X_{ \xbb },$ $ \xbb < \xbk.$ Thus, from $X'_{ \xbk },$ we
cannot go anywhere interesting (recall that the intersections are not
candidates
for a $Y$ giving a contradiction).

Consider now $X_{ \xba }.$ We can go up to any $X_{ \xba +n},$ but not to
any $X_{ \xbl },$ $ \xba < \xbl,$ as
$d_{ \xbl }$ is missing, neither to $X'_{ \xbk },$ as a is missing. And we
will be stopped by
the first $ \xbl > \xba,$ as $x_{ \xbl }$ will be missing to go beyond
$X_{ \xbl }.$ Analogous observations
hold for the remaining intersections (1), (6.1), (7.2). But in all these
sets we
can reach, we will not destroy minimality of any element of $X_{ \xba }$
(or of the
intersections).

Consequently, the only candidates for failure will all start with $U.$ As
the only
element of $U$ not globally minimal is $c,$ such failure has to have $c
\xbe Y- \xbm (Y),$ so
$Y$ has to be $X'_{ \xbk }.$ Suppose we omit one of the $X_{ \xba }$ in
the sequence going up to
$X'_{ \xbk }.$ If $ \xbk \xcg \xbl > \xba,$ we cannot reach $X_{ \xbl }$
and beyond, as $x'_{ \xba }$ will be missing.
But we cannot go to $X_{ \xba +n}$ either, as $x_{ \xba +1}$ is missing.
So we will be stopped
at $X_{ \xba }.$ Thus, to see failure, we need the full sequence
$U=X_{0},$ $X'_{ \xbk }=Y_{ \xbk },$
$Y_{ \xba }=X_{ \xba }$ for $0< \xba < \xbk.$

\paragraph{
$( \xbm Cum \xbk )$ fails:
}

$\hspace{0.01em}$


\label{Section ( mCum k) fails:}

The full sequence $U=X_{0},$ $X'_{ \xbk }=Y_{ \xbk },$ $Y_{ \xba }=X_{
\xba }$ for $0< \xba < \xbk $ shows this, as
$c \xbe \xbm (U) \xcs X'_{ \xbk },$ but $c \xce \xbm (X'_{ \xbk }).$

Consequently, the example satisfies $( \xcS ),$ $( \xbm CUM),$ $( \xbm
Cumt \xba )$ for $ \xba < \xbk,$ and
$( \xbm Cum \xbk )$ fails.

$ \xcz $
\\[3ex]

To put our work more into perspective, we mention and prove now
some positive results about the $( \xbm Cum \xba )$ and $( \xbm Cumt \xba
).$
\index{Fact Cum-Alpha}

\bfa

$\hspace{0.01em}$


\label{Fact Cum-Alpha}

We summarize some properties of $( \xbm Cum \xba )$ and $( \xbm Cumt \xba
)$ - sometimes with some
redundancy. Unless said otherwise, $ \xba,$ $ \xbb $ etc. will be
arbitrary ordinals.

For (1) to (6) $( \xbm PR)$ and $( \xbm \xcc )$ are assumed to hold, for
(7) only
$( \xbm \xcc ).$

(1) Downward:

(1.1) $( \xbm Cum \xba )$ $ \xcp $ $( \xbm Cum \xbb )$ for all $ \xbb \xck
\xba $

(1.2) $( \xbm Cumt \xba )$ $ \xcp $ $( \xbm Cumt \xbb )$ for all $ \xbb
\xck \xba $

\efa

(2) Validity of $( \xbm Cum \xba )$ and $( \xbm Cumt \xba )$:

(2.1) All $( \xbm Cum \xba )$ hold in smooth preferential structures

(2.2) All $( \xbm Cumt \xba )$ hold in transitive smooth preferential
structures

(2.3) $( \xbm Cumt \xba )$ for $0< \xba $ do not necessarily hold in
smooth structures without
transitivity, even in the presence of $( \xcS )$

(3) Upward:

(3.1) $( \xbm Cum \xbb )$ $+$ $( \xcv )$ $ \xcp $ $( \xbm Cum \xba )$ for
all $ \xbb \xck \xba $

(3.2) $( \xbm Cumt \xbb )$ $+$ $( \xcv )$ $ \xcp $ $( \xbm Cumt \xba )$
for all $ \xbb \xck \xba $

(3.3) $\{( \xbm Cumt \xbb ): \xbb < \xba \}$ $+$ $( \xbm CUM)$ $+$ $( \xcS
)$ $ \xcP $ $( \xbm Cum \xba )$ for $ \xba >0.$

(4) Connection $( \xbm Cum \xba )/( \xbm Cumt \xba )$:

(4.1) $( \xbm Cumt \xba )$ $ \xcp $ $( \xbm Cum \xba )$

(4.2) $( \xbm Cum \xba )$ $+$ $( \xcS )$ $ \xcP $ $( \xbm Cumt \xba )$

(4.3) $( \xbm Cum \xba )$ $+$ $( \xcv )$ $ \xcp $ $( \xbm Cumt \xba )$

(5) $( \xbm CUM)$ and $( \xbm Cumi)$:

(5.1) $( \xbm CUM)$ $+$ $( \xcv )$ entail:

(5.1.1) $ \xbm (A) \xcc B$ $ \xcp $ $ \xbm (A \xcv B)= \xbm (B)$

(5.1.2) $ \xbm (X) \xcc U,$ $U \xcc Y$ $ \xcp $ $ \xbm (Y \xcv X)= \xbm
(Y)$

(5.1.3) $ \xbm (X) \xcc U,$ $U \xcc Y$ $ \xcp $ $ \xbm (Y) \xcs X \xcc
\xbm (U)$

(5.2) $( \xbm Cum \xba )$ $ \xcp $ $( \xbm CUM)$ for all $ \xba $

(5.3) $( \xbm CUM)$ $+$ $( \xcv )$ $ \xcp $ $( \xbm Cum \xba )$ for all $
\xba $

(5.4) $( \xbm CUM)$ $+$ $( \xcs )$ $ \xcp $ $( \xbm Cum0)$

(6) $( \xbm CUM)$ and $( \xbm Cumt \xba )$:

(6.1) $( \xbm Cumt \xba )$ $ \xcp $ $( \xbm CUM)$ for all $ \xba $

(6.2) $( \xbm CUM)$ $+$ $( \xcv )$ $ \xcp $ $( \xbm Cumt \xba )$ for all $
\xba $

(6.3) $( \xbm CUM)$ $ \xcP $ $( \xbm Cumt \xba )$ for all $ \xba >0$

(7) $( \xbm Cum0)$ $ \xcp $ $( \xbm PR)$
\index{Fact Cum-Alpha Proof}

\paragraph{
Proof of Fact Cum-Alpha
}

$\hspace{0.01em}$


\label{Section Proof of Fact Cum-Alpha}

We prove these facts in a different order: (1), (2), (5.1), (5.2), (4.1),
(6.1),
(6.2), (5.3), (3.1), (3.2), (4.2), (4.3), (5.4), (3.3), (6.3), (7).

(1.1)

For $ \xbb < \xbg \xck \xba $ set $X_{ \xbg }:=X_{ \xbb }.$ Let the
prerequisites of $( \xbm Cum \xbb )$ hold. Then for
$ \xbg $ with $ \xbb < \xbg \xck \xba $ $ \xbm (X_{ \xbg }) \xcc X_{ \xbb
}$ by $( \xbm \xcc ),$ so the prerequisites
of $( \xbm Cum \xba )$ hold, too, so by $( \xbm Cum \xba )$ $ \xcS \{X_{
\xbd }: \xbd \xck \xbb \} \xcs \xbm (U)$ $=$
$ \xcS \{X_{ \xbd }: \xbd \xck \xba \} \xcs \xbm (U)$ $ \xcc $ $ \xbm (X_{
\xba })$ $=$ $ \xbm (X_{ \xbb }).$

(1.2)

Analogous.

(2.1)

Proof by induction.

$( \xbm Cum0)$ Let $ \xbm (X_{0}) \xcc U,$ suppose there is $x \xbe \xbm
(U) \xcs (X_{0}- \xbm (X_{0})).$ By smoothness,
there is $y \xeb x,$ $y \xbe \xbm (X_{0}) \xcc U,$ $contradiction$ (The
same arguments works for copies: all
copies of $x$ must be minimized by some $y \xbe \xbm (X_{0}),$ but at
least one copy of $x$
has to be minimal in $U.)$

Suppose $( \xbm Cum \xbb )$ hold for all $ \xbb < \xba.$ We show $( \xbm
Cum \xba ).$ Let the prerequisites
of $( \xbm Cum \xba )$ hold, then those for $( \xbm Cum \xbb ),$ $ \xbb <
\xba $ hold, too. Suppose there is
$x \xbe \xbm (U) \xcs \xcS \{X_{ \xbg }: \xbg \xck \xba \}- \xbm (X_{ \xba
}).$ So by $( \xbm Cum \xbb )$ for $ \xbb < \xba $ $x \xbe \xbm (X_{ \xbb
})$
moreover $x \xbe \xbm (U).$ By smoothness, there is $y \xbe \xbm (X_{ \xba
}) \xcc U \xcv \xcV \{X_{ \xbb }: \xbb < \xba \},$ $y \xeb x,$
but this is a contradiction. The same argument works again for copies.

(2.2)

We use the following Fact:
Let, in a smooth transitive structure, $ \xbm (X_{ \xbb })$ $ \xcc $ $U
\xcv \xcV \{X_{ \xbg }: \xbg < \xbb \}$ for all $ \xbb \xck \xba,$
and let $x \xbe \xbm (U).$ Then there is no $y \xeb x,$ $y \xbe U \xcv
\xcV \{X_{ \xbg }: \xbg \xck \xba \}.$

Proof of the Fact by induction:
$ \xba =0:$ $y \xbe U$ is impossible: if $y \xbe X_{0},$ then if $y \xbe
\xbm (X_{0}) \xcc U,$ which is impossible,
or there is $z \xbe \xbm (X_{0}),$ $z \xeb y,$ so $z \xeb x$ by
transitivity, but $ \xbm (X_{0}) \xcc U.$
Let the result hold for all $ \xbb < \xba,$ but fail for $ \xba,$
so $ \xCN \xcE y \xeb x.y \xbe U \xcv \xcV \{X_{ \xbg }: \xbg < \xba \},$
but $ \xcE y \xeb x.y \xbe U \xcv \xcV \{X_{ \xbg }: \xbg \xck \xba \},$
so $y \xbe X_{ \xba }.$
If $y \xbe \xbm (X_{ \xba }),$ then $y \xbe U \xcv \xcV \{X_{ \xbg }: \xbg
< \xba \},$ but this is impossible, so
$y \xbe X_{ \xba }- \xbm (X_{ \xba }),$ let by smoothness $z \xeb y,$ $z
\xbe \xbm (X_{ \xba }),$ so by transitivity $z \xeb x,$ $contradiction.$
The result is easily modified for the case with copies.

Let the prerequisites of $( \xbm Cumt \xba )$ hold, then those of the Fact
will hold,
too. Let now $x \xbe \xbm (U) \xcs (X_{ \xba }- \xbm (X_{ \xba })),$ by
smoothness, there must be $y \xeb x,$
$y \xbe \xbm (X_{ \xba }) \xcc U \xcv \xcV \{X_{ \xbg }: \xbg < \xba \},$
contradicting the Fact.

(2.3)

Let $ \xba >0,$ and consider the following structure over $\{a,b,c\}:$
$U:=\{a,c\},$
$X_{0}:=\{b,c\},$ $X_{ \xba }:= \Xl:=X_{1}:=\{a,b\},$ and their
intersections, $\{a\},$ $\{b\},$ $\{c\},$ $ \xCQ $ with
the order $c \xeb b \xeb a$ (without transitivity). This is preferential,
so $( \xbm PR)$ and
$( \xbm \xcc )$ hold.
The structure is smooth for $U,$ all $X_{ \xbb },$ and their
intersections.
We have $ \xbm (X_{0}) \xcc U,$ $ \xbm (X_{ \xbb }) \xcc U \xcv X_{0}$ for
all $ \xbb \xck \xba,$ so $ \xbm (X_{ \xbb }) \xcc U \xcv \xcV \{X_{ \xbg
}: \xbg < \xbb \}$
for all $ \xbb \xck \xba $ but $X_{ \xba } \xcs \xbm (U)=\{a\} \xcC \{b\}=
\xbm (X_{ \xba })$ for $ \xba >0.$

(5.1)

(5.1.1) $ \xbm (A) \xcc B$ $ \xcp $ $ \xbm (A \xcv B) \xcc \xbm (A) \xcv
\xbm (B) \xcc B$ $ \xcp_{( \xbm CUM)}$ $ \xbm (B)= \xbm (A \xcv B).$

(5.1.2) $ \xbm (X) \xcc U \xcc Y$ $ \xcp $ (by (1)) $ \xbm (Y \xcv X)=
\xbm (Y).$

(5.1.3) $ \xbm (Y) \xcs X$ $=$ (by (2)) $ \xbm (Y \xcv X) \xcs X$ $ \xcc $
$ \xbm (Y \xcv X) \xcs (X \xcv U)$ $ \xcc $ (by $( \xbm PR))$
$ \xbm (X \xcv U)$ $=$ (by (1)) $ \xbm (U).$

(5.2)

Using (1.1), it suffices to show $( \xbm Cum0)$ $ \xcp $ $( \xbm CUM).$
Let $ \xbm (X) \xcc U \xcc X.$ By $( \xbm Cum0)$ $X \xcs \xbm (U) \xcc
\xbm (X),$ so by $ \xbm (U) \xcc U \xcc X$ $ \xcp $ $ \xbm (U) \xcc \xbm
(X).$
$U \xcc X$ $ \xcp $ $ \xbm (X) \xcs U \xcc \xbm (U),$ but also $ \xbm (X)
\xcc U,$ so $ \xbm (X) \xcc \xbm (U).$

(4.1)

Trivial.

(6.1)

Follows from (4.1) and (5.2).

(6.2)

Let the prerequisites of $( \xbm Cumt \xba )$ hold.

We first show by induction $ \xbm (X_{ \xba } \xcv U) \xcc \xbm (U).$

Proof:

$ \xba =0:$ $ \xbm (X_{0}) \xcc U$ $ \xcp $ $ \xbm (X_{0} \xcv U)= \xbm
(U)$ by (5.1.1).
Let for all $ \xbb < \xba $ $ \xbm (X_{ \xbb } \xcv U) \xcc \xbm (U) \xcc
U.$ By prerequisite,
$ \xbm (X_{ \xba }) \xcc U \xcv \xcV \{X_{ \xbb }: \xbb < \xba \},$ thus $
\xbm (X_{ \xba } \xcv U)$ $ \xcc $ $ \xbm (X_{ \xba }) \xcv \xbm (U)$ $
\xcc $ $ \xcV \{U \xcv X_{ \xbb }: \xbb < \xba \},$

so $ \xcA \xbb < \xba $ $ \xbm (X_{ \xba } \xcv U) \xcs (U \xcv X_{ \xbb
})$ $ \xcc $ $ \xbm (U)$ by (5.1.3), thus $ \xbm (X_{ \xba } \xcv U) \xcc
\xbm (U).$

Consequently, under the above prerequisites, we have $ \xbm (X_{ \xba }
\xcv U)$ $ \xcc $ $ \xbm (U)$ $ \xcc $
$U$ $ \xcc $ $U \xcv X_{ \xba },$ so by $( \xbm CUM)$ $ \xbm (U)= \xbm
(X_{ \xba } \xcv U),$ and, finally,
$ \xbm (U) \xcs X_{ \xba }= \xbm (X_{ \xba } \xcv U) \xcs X_{ \xba } \xcc
\xbm (X_{ \xba })$ by $( \xbm PR).$

Note that finite unions take us over the limit step, essentially, as all
steps collaps, and $ \xbm (X_{ \xba } \xcv U)$ will always be $ \xbm (U),$
so there are no real
changes.

(5.3)

Follows from (6.2) and (4.1).

(3.1)

Follows from (5.2) and (5.3).

(3.2)

Follows from (6.1) and (6.2).

(4.2)

Follows from (2.3) and (2.1).

(4.3)

Follows from (5.2) and (6.2).

(5.4)

$ \xbm (X) \xcc U$ $ \xcp $ $ \xbm (X) \xcc U \xcs X \xcc X$ $ \xcp $ $
\xbm (X \xcs U)= \xbm (X)$ $ \xcp $
$X \xcs \xbm (U)=(X \xcs U) \xcs \xbm (U) \xcc \xbm (X \xcs U)= \xbm (X)$

(3.3)

See Example \ref{Example Inf-Cum-Alpha}.

(6.3)

This is a consequence of (3.3).

(7)

Trivial. Let $X \xcc Y,$ so by $( \xbm \xcc )$ $ \xbm (X) \xcc X \xcc Y,$
so by $( \xbm Cum0)$ $X \xcs \xbm (Y) \xcc \xbm (X).$

$ \xcz $
\\[3ex]

\end{document}